\documentclass[12pt,draftcls,onecolumn]{IEEEtran}

\usepackage{amsfonts}
\usepackage{amssymb}
\usepackage{amsmath}
\usepackage{mathptmx}
\usepackage{cancel}
\usepackage{pstricks}
\usepackage{psfrag}
\usepackage{mathrsfs}
\usepackage{graphicx}
\usepackage{pgf,tikz}
\usepackage{algorithm}
\usepackage{algpseudocode}
\usepackage{soul}
\def\be{\begin{equation}}
\def\ee{\end{equation}}
\def\bea{\begin{eqnarray}}
\def\eea{\end{eqnarray}}
\def\beann{\begin{eqnarray*}}
\def\eeann{\end{eqnarray*}}

\newcommand{\rank}{{\rm rank}}

\newcommand{\real}{{\mathbb{R}}}

\def\spacingset#1{\def\baselinestretch{#1}\small\normalsize}
\spacingset{1.35}

\newtheorem{lemma}{Lemma}
\newtheorem{theorem}{Theorem}
\newtheorem{remark}{Remark}

\newtheorem{example}{Example}[section]

\def\be{\begin{equation}}
\def\ee{\end{equation}}
\def\bea{\begin{eqnarray}}
\def\eea{\end{eqnarray}}
\def\beann{\begin{eqnarray*}}
\def\eeann{\end{eqnarray*}}

\def\proof{\noindent{\bf{\em Proof:}\ \ }}
\def\QED{\mbox{\rule[0pt]{1.5ex}{1.5ex}}}
\def\endproof{\hspace*{\fill}~\QED\par\endtrivlist\unskip}
\def\endex{\hspace*{\fill}~\Box\par\endtrivlist\unskip}

\def\second{{\prime \prime}}
\def\tilda{{\!\!\!\!\phantom{P}^\thicksim}}

\newcommand{\ima}{\operatorname{im}}
\newcommand{\normrank}{\operatorname{normrank}}
\newcommand{\diag}{\operatorname{diag}}
\newcommand{\defi}{\stackrel{\text{\tiny def}}{=}}

\newcommand{\complex}{{\mathbb{C}}}

\def\gR{{\cal R}}

\def\gV{{\cal V}}

\def\gX{{\cal X}}

\def\gZ{{\cal Z}}

\def\bmat{\left[ \begin{array}}
\def\emat{\end{array} \right]}
\def\endex{\hspace*{\fill}~$\square$\par\endtrivlist\unskip}
\def\bmat{\left[ \begin{array}}
\def\emat{\end{array} \right]}
\def\bsmat{\left[ \begin{smallmatrix}}
\def\esmat{\end{smallmatrix} \right]}

\def\gR{{\cal R}}

\def\gV{{\cal V}}

\def\gX{{\cal X}}
\def\i{{i}}
\newcommand{\spanR}{\operatorname{span}}

\begin{document}
\title{\LARGE{The geometry of the generalized  algebraic Riccati equation and of the singular Hamiltonian system}}

\author{Lorenzo~Ntogramatzidis, and Augusto Ferrante
 % <-this % stops r space
 \thanks{L. Ntogramatzidis is with the Department of Mathematics and
Statistics, Curtin University, Perth,
Australia. E-mail: {\tt L.Ntogramatzidis@curtin.edu.au. }}

% \thanks{Jean-Fran\c{c}ois Tr{\'e}gou{\"e}t is with Universit\'{e} de Lyon, Laboratoire Amp\`{e}re CNRS UMR 5005, INSA-Lyon; F-69621, Villeurbanne, France. E-mail: {\tt jean-francois.tregouet@insa-lyon.fr. } (Research partially carried out at Curtin University).}

%\thanks{R. Schmid is with the Department of Electrical and Electronic Engineering,
%The University of Melbourne, Parkville, VIC 3010, Australia. E-mail: {\tt rschmid@unimelb.edu.au. }}

 \thanks{A. Ferrante is with the Dipartimento di Ingegneria dell'Informazione, Universit\`a di Padova, 
via Gradenigo, 6/B -- I-35131 Padova, Italy. E-mail: {\tt augusto@dei.unipd.it. }}

% <-this % stops a space
%\thanks{Manuscript received March 23, 2010; revised XXX}}
}

% The paper headers
%\markboth{DRAFT}{Shell \MakeLowercase{\textit{et al.}}: Bare Demo of IEEEtran.cls for Journals}

\maketitle

\vspace{-1cm}

\IEEEpeerreviewmaketitle

\begin{abstract}
This paper analyzes the properties of the solutions of the generalized continuous algebraic Riccati equation from a geometric perspective.  This analysis reveals the presence of a subspace that 
may provide an appropriate degree of freedom to
stabilize the system in the related optimal control problem even in cases
where the Riccati equation does not admit a stabilizing
solution.
\end{abstract}

%\newpage
%\tableofcontents
%\newpage

\section{Introduction}
\label{secintro} 
This paper investigates the geometric properties of the set of solutions of the so-called {\em constrained generalized continuous algebraic Riccati equation} associated with the infinite-horizon linear quadratic (LQ) optimal control problem, when the matrix $R$  weighting the input in the cost function is allowed to be singular. 
 This problem, often referred to as the singular LQ problem, has a long history. It has been investigated in several papers and with different techniques, see
\cite{Hautus-S-83,Willems-KS-86,Saberi-S-87,Prattichizzo-MN-04,Kalaimani-BC-13} and the references therein.
See also the monographs \cite{Abou-Kandil-FIJ-03,Lancaster-95,Ionescu-OW-99,Saberi-SC-95} for a more general discussion.

 In particular, in the foundational contributions
\cite{Hautus-S-83} and \cite{Willems-KS-86} it was proved that an optimal solution of the singular LQ problem exists for all initial conditions if the class of controls is extended to include
distributions.  A different perspective was offered in \cite{Prattichizzo-MN-04}, where a geometric approach was employed on the Hamiltonian differential equation to study the subspace of initial conditions for which the control law is impulse-free.

In the discrete time this issue does not arise, and  it is an established fact that  the solution of regular and singular infinite-horizon LQ problem can be found resorting to the so-called  {\em constrained generalized discrete algebraic Riccati equation}, see  \cite{Ferrante-N-12}. 
Considerable effort has been devoted  | also in recent years |  in providing a geometric characterization of the set of solutions of this discrete Riccati equation, see e.g. \cite{Stoorvogel-S-98} and \cite{Ferrante-N-12}.  A similar characterization for the continuous-time generalized Riccati equation 
has never been considered.

 There are several reasons for considering this equation and for analyzing the geometric structure of its solutions. The first, which is our main motivation, is given by the recent results connecting the continuous time generalized Riccati equation with LQ optimal control problems \cite{Ferrante-N-14}.
Another reason derives from the fact that this equation is a particular case of an even more general type of Riccati equation that arises in the literature that flourished in the past twenty years on stochastic optimal control, see e.g. \cite{Abou-Kandil-FIJ-03,Damm-04,Damm-H-01,Dragan-MS-10,Freiling-H-03,Freiling-H-04} and the references cited therein as well as \cite{zorzi1,zorzi2,zorzi3} for  the dual version in filtering problems. These research lines may benefit of our contribution. In fact, the natural approach in this field is based on the study of the corresponding Hamiltonian system, so that our new geometric results may furnish a powerful point of view to deal with these problems and with the associated numerical analysis.

%Despite the lack of a real understanding of the links with control problems, 
In 
 \cite{Ionescu-O-96-1} the constrained generalized continuous algebraic Riccati equation was defined, in analogy with the discrete case, by replacing  the inverse of the matrix $R$ appearing in the standard Riccati equation with its pseudo-inverse. 
In particular, this paper offers a characterization in terms of deflating subspaces of the Hamiltonian pencil of the conditions under which the constrained generalized Riccati equation has a stabilizing solution.

To our best knowledge, the recent papers \cite{Ferrante-N-14,Ferrante-N-14-1} were the first attempts to link this equation to singular LQ optimal control problems. 
In \cite{Ferrante-N-14,Ferrante-N-14-1} it was shown that the existence of symmetric solutions of the constrained generalized continuous-time Riccati equation is equivalent to the existence of impulse-free solutions of the associated singular LQ problem from any initial condition. This means, in particular, that an optimal control can always be expressed as a { state-feedback. 
Now that the connection between the constrained generalized continuous-time algebraic Riccati equation and the singular LQ problem has been explained, the important issue arises of analyzing the set of solutions of such equation and the relations of each such solution with the corresponding LQ control problem.}
 
 In this paper a geometric analysis is carried out on the structure of the symmetric solutions of the constrained generalized continuous-time algebraic Riccati equation. This analysis leads to the following main contributions.
First, we show that the dynamics of the closed-loop system can be divided into a part that depends on the particular solution $X$ that we are considering,   and one 
which is independent of it. We also show that the latter dynamics, which is not necessarily stable, is confined to an output nulling subspace, so that it does not contribute to the cost function. The spectrum associated with the reachable part of this dynamics can therefore be assigned without affecting the optimality of the cost. As a consequence, we show that the LQ optimal control problem may admit a stabilizing solution even in cases in which the generalized continuous-time Riccati equation does not admit a stabilizing solution. This is a new feature that has no parallel in the regular LQ problems.
We finally address the analysis of the structure of the corresponding Hamiltonian system and its relations with the generalized  algebraic Riccati equations and the singular LQ optimal control problems: we show that differently from the regular case,  only  the eigenvalues of the closed-loop dynamics that depend on the particular solution $X$ correspond -- together with their mirrored values -- to the invariant zeros of the Hamiltonian system. 
An anonymous reviewer has pointed out that 
some of the results of this paper  may be  alternatively obtained
by performing a preliminary transformation that brings the system in  the so-called {\em special coordinate basis} of \cite{Saberi-S-87}. We believe that a direct derivation of these results will provide additional insight to some readers  as it connects the results with the structure of the Hamiltonian system.

 \section{The generalized Riccati equation and Linear Quadratic optimal control}
\label{LQ}
Let $Q\in \real^{n \times n}$, $S \in \real^{n \times m}$, $R \in \real^{m \times m}$.  We make the following standing assumption:
 \begin{equation}
 \label{equno}
\Pi \defi \bmat{cc}  Q & S \\[-1mm] S^\top & R \emat=\Pi^\top \ge 0.
 \end{equation}
The triplet $\Sigma \defi (A,B,\Pi)$ is referred to as {\em Popov} triple.

From the properties of the Schur complement, we recall that the condition $\Pi=\Pi^\top \ge 0$ is equivalent to the simultaneous satisfaction of the three conditions
\begin{itemize}
\item $R\ge 0$;
\item $\ker R \subseteq \ker S$;
\item $Q-S\,R^\dagger S^\top \ge 0$;
\end{itemize}
Dually, $\Pi\ge 0$ if and only if
\begin{itemize}
\item $Q\ge 0$;
\item $\ker Q \subseteq \ker S^\top$;
\item $R-S^\top\,Q^\dagger S \ge 0$.
\end{itemize}
See e.g. \cite{Rami-CZ-02} or \cite{Ferrante-N-12}  for a proof. 
From these considerations it follows also that if $\Pi=\Pi^\top \ge 0$, then $S\,R^\dagger\,R=S$ and $S^\top\,Q^\dagger\,Q=S^\top$.

The classic LQ problem can be stated as the problem of finding the control $u(t)$, $t \ge 0$, that minimizes 
   \begin{equation}
 \label{costinf}
J_\infty(x_{ 0},u)=\int_0^\infty [\begin{array}{cc} x^\top(t) & u^\top(t) \end{array} ] \bmat{cc} Q & S \\[-1mm] S^\top & R \emat \bmat{c} x(t) \\[-1mm] u(t) \emat\,dt
\end{equation}
subject to the constraint 
 \begin{equation}
 \label{eqsys}
 \dot{x}(t)=A\,x(t)+B\,u(t), \qquad x(0)=x_{ 0} \in \real^n
 \end{equation}
 where $A\in \real^{n \times n}$ and $B \in \real^{n \times m}$. 
When $R$ is positive definite, the optimal control (when it exists) does not include distributions, since in such a case an impulsive control $u$ will always cause $J_\infty(x_{ 0},u)$ to be unbounded for any $x_{ 0}\in \real^n$. If $R$ is only positive semidefinite, in general the optimal solution can contain the Dirac delta distribution and its derivatives.  In the very recent literature, it has been shown that important links exist between the existence  of the solutions of the constrained generalized continuous algebraic Riccati equation (often denoted by the acronym CGCARE and formally introduced in the next section) and the non-impulsive optimal solutions of the infinite-horizon LQ problem, \cite{Ferrante-N-14,Ferrante-N-14-1}. This point represents a crucial difference between the discrete and the continuous time. Indeed, while in the discrete time the existence of symmetric positive semidefinite solutions of the constrained generalized discrete algebraic Riccati equation is equivalent to the solvability of the infinite-horizon LQ problem, in the continuous-time case this correspondence holds for the so-called {\em regular} solutions, i.e., the optimal controls of the LQ problem that do not contain distributions.

LQ problems have been found to be very important as control problems in their own right. On the other hand, in the last thirty years the LQ problem has been often used as a building block to solve different, and usually more articulated, optimal control problems. For example, in the so-called  $H_2$ problem  \cite{Stoorvogel-92_2}  the index to be minimized is the norm of the output of the system
\[
y(t)=C\,x(t)+D\,u(t).
\]
The corresponding LQ problem is obtained by defining $Q=C^\top\,C$, $S=C^\top\,D$ and $R=D^\top\,D$. Since the very vast majority of systems (for example virtually all mechanical systems) are strictly proper, then the corresponding LQ problem is usually singular.

\section{Generalized CARE}
Consider the following matrix equation\footnote{The symbol $M^\dagger$ denotes the Moore-Penrose pseudo-inverse of matrix $M$.}
 \begin{equation}
 X\,A+A^\top\,X-(S+X\,B)\,R^{\dagger}\,(S^\top+B^\top X)+Q=0, \label{gcare}
 \end{equation}
where the matrices $Q, A\in \real^{n \times n}$, $B,S \in \real^{n \times m}$, $R \in \real^{m \times m}$ are as defined in the previous section. Equation (\ref{gcare}), where $R$ is allowed to be singular, is often referred to as the {\em generalized continuous algebraic Riccati equation} GCARE($\Sigma$). 
Equation (\ref{gcare}) along with the condition
 \begin{equation}
 \ker R \subseteq \ker (S+X\,B), \label{kercond}
 \end{equation}
 will be referred to as {\em constrained generalized continuous algebraic Riccati equation}, and denoted by CGCARE($\Sigma$). In view of the positive semidefiniteness of $\Pi$, as already observed in Section \ref{LQ}, we have $\ker R \subseteq \ker S$, which implies that (\ref{kercond}) is equivalent to 
$\ker R \subseteq \ker (X\,B)$. 
 The following notation is used throughout the paper. First, let $G \defi I_m-R^\dagger R$ be  the orthogonal projector that projects onto $\ker R$.  Moreover, we consider a non-singular matrix $T=[T_{1}\mid T _{2}]$ where
$\ima T_{1}=\ima R$ and $\ima T _{2}=\ima G$, and we define $B_{1}\defi BT_{1}$ and $B _{2} \defi BT _{2}$. Finally, to any $X=X^\top \in \real^{n \times n}$  we associate the following matrices
\bea
Q _{X}& \defi & Q+A^\top X+X\,A, \qquad
S _{X}   \defi   S+X\, B, \label{defgx} \\
K _{X} & \defi & R^\dagger\, S _{X}^\top, \qquad A _{X} \defi  A-B\,K _{X}, \qquad 
  \Pi _{X}  \defi  \left[ \begin{array}{cc} Q _{X} & S _{X} \\[-1mm] S _{X}^\top & R \end{array} \right].
 \eea 
When $X$ is a solution of CGDARE($\Sigma$), then $K _{X}$ is the corresponding gain matrix, and $A _{X}$ the associated closed-loop matrix. 

\begin{remark}
{
A symmetric and positive semidefinite solution of the generalized discrete-time algebraic Riccati equation also solves the constrained generalized discrete-time algebraic Riccati equation.
This fact does not hold in the continuous time, i.e.,  not all symmetric and positive semidefinite solutions of GCARE($\Sigma$) are also solutions of CGCARE($\Sigma$). 
}
\end{remark}
\section{Characterization of the solutions of CGCARE}
The purpose of this section is to provide a geometric characterization for the set of solutions of the generalized continuous algebraic Riccati equation.
To this end, we first recall some concepts of classical geometric control theory that will be used in the sequel. More details can be found e.g. in \cite{Trentelman-SH-01}. Consider a system described by (\ref{eqsys}) along with the output equation $y(t) = C\,x(t)+D\,u(t)$, that we concisely identify with the quadruple $\Sigma_0=(A,B,C,D)$.

 The {\em invariant zeros} of $\Sigma_0$, here denoted by $\gZ({A}, {B}, {C}, {D})$, are the values $s \in \complex$ such that the rank of the Rosenbrock system matrix pencil $\bsmat A-s\,I_n & B \\[1mm] C & D \esmat$ is smaller than its normal rank, \cite[Def. 3.16]{Zhou-DG-96}. We recall that the {\em reachable subspace} is 
 $\gR_0=\ima [\begin{array}{ccccccccc} B && A\,B && \ldots && A^{n-1}\,B \end{array}]$, and
coincides with the smallest $A$-invariant subspace of $\real^n$ containing the image of $B$, i.e. $\gR_0=\langle A\,|\,\ima B \rangle$. 
An {\em output-nulling subspace} $\gV$ of $\Sigma_0$ is a subspace of $\real^n$ %which satisfies the inclusion
%\bea
%\label{def}
 %\bmat{c} A \\ C \emat \,\gV \subseteq (\gV \oplus \{0\}) + \ima\, \bmat{c} B \\ D \emat,
%\eea
%which is equivalent to the existence of 
for which there exists a matrix $F\,{\in}\,\mathbb{R}^{m\,{\times}\,n}$ such that
$(A+B\,F)\,\gV\subseteq \gV \subseteq \ker (C+D\,F)$. Any real matrix $F$ satisfying these inclusions is referred to as a {\it friend \/} of $\gV$. We denote by $\mathfrak{F}(\gV)$ the set of friends of $\gV$. 
We denote by $\gV^\star$ the largest output-nulling subspace of $\Sigma_0$, which represents the set of all initial states $x_0$ of $\Sigma_0$ for which a control input exists such that the corresponding output function is identically zero. Such an input function can always be implemented as a static state feedback of the form $u(t)=F\,x(t)$ where $F \in \mathfrak{F}(\gV^\star)$.
 The so-called {\em output-nulling reachability subspace} on $\gV^\star$, herein denoted with $\gR^\star$, is the smallest $(A\,{+}\,B\,F)$-invariant subspace of $\real^n$ containing the subspace $\gV^\star\,{\cap}\,B\,\ker\,D$, where $F\,{\in}\,\mathfrak{F}(\gV^\star)$, i.e.,
$\gR^\star=\langle A+B\,F\,|\, \gV^\star \cap B\,\ker D\rangle$ where $F \in \mathfrak{F}(\gV^\star)$. %
Let $F\in \mathfrak{F}(\gV^\star)$. The closed-loop spectrum (viewed as a multiset, with aggregation denoted by $\uplus$) can be partitioned as
$\sigma(A+B\,F)=\sigma(A+B\,F\,|\,\gV^\star)\uplus \sigma(A+B\,F\,|\,\gX/\gV^\star)$, where $\sigma(A+B\,F\,|\,\gV^\star)$ is the spectrum of $A+B\,F$ restricted to $\gV^\star$ and $\sigma(A+B\,F\,|\,\gX/\gV^\star)$ is the spectrum of the mapping induced by $A+B\,F$ on the quotient space $\gX/\gV^\star$. 
 The eigenvalues of $A+B\,F$ restricted to $\gV^\star$ can be further split into two disjoint sets: the eigenvalues of $\sigma(A+B\,F |\gR^\star)$ are all freely assignable
 with a suitable choice of $F$ in $\mathfrak{F}(\gV^\star)$. The eigenvalues in $\sigma\,(A+B\,F | {\gV^\star}/{\gR^\star})$ -- which coincide with the {\em invariant zeros} of $\Sigma_0$, see e.g. \cite[Theorem 7.19]{Trentelman-SH-01}  
  -- are fixed for all the choices of $F$ in $\mathfrak{F}(\gV^\star)$.

Since $\Pi$ is assumed symmetric and positive semidefinite, we can consider a factorization of the form
 \begin{equation}
\label{pifact}
\Pi=\left[ \begin{array}{cc} Q & S \\[-1mm] S^\top & R\end{array} \right]=\left[ \begin{array}{cc} C^\top \\[-1mm] D^\top\end{array} \right][ \begin{array}{cc} C & D \end{array} ],
 \end{equation}
where $Q=C^\top C$, $S=C^\top D$ and $R=D^\top D$.  Let us define $G(s)  \defi  C\,(s\,I_n-A)^{-1}B+D$. Let $G^\tilda(s) \defi G^\top(-s)$. The ``spectrum'' or ``spectral density'' $\Phi(s) \defi G^\tilda(s)\,G(s)$ %associated with $\Sigma$ 
can be written as
\beann
\Phi(s)=[ \begin{array}{cc} B^\top (-s\,I_n-A^\top)^{-1} & I_m \end{array} ]\,
\left[ \begin{array}{cc} Q & S \\[-1mm] S^\top & R \end{array} \right]\,\left[ \begin{array}{cc} (s\,I_n-A)^{-1}\,B \\[-1mm] I_m \end{array} \right],
\eeann
which is also referred to as {\em Popov function}. We recall the following classical result. % associated with GCARE($\Sigma$). 
%Recall that we have defined $\Pi _{X}=\bsmat Q _{X} & S _{X} \\[1mm] S _{X}^\top & R \esmat= \bsmat Q+A^\top X+X\,A & S+X\,B \\[1mm] S^\top+B^\top X & R\esmat$. 
%Let us also define
%\beann
%L(X) \defi \Pi _{X}-\Pi=\left[ \begin{array}{cc} A^\top X+X\,A & X\,B \\ B^\top X & 0 \end{array} \right].
%\eeann
%Notice that $L(X)$ is a linear function of $X$.
\begin{lemma}%{{\bf (\cite[p.322]{Stoorvogel-S-98}, see e.g. \cite{Colaneri-F-SCL} for a detailed proof).}
For any $X=X^\top \in\real^{n\times n}$, there holds 
\begin{equation}
\label{alpha}
\Phi(s)=[ \begin{array}{cc} B^\top (-s\,I_n-A^\top)^{-1} & I_m \end{array} ]\,
\Pi _{X}\,\left[ \begin{array}{cc} (s\,I_n-A)^{-1}\,B \\[-1mm] I_m \end{array} \right].
\end{equation}
\end{lemma}
%The proof of this result is standard, see e.g. \cite[p.322]{Stoorvogel-S-98}.
%
\proof 
The statement follows on noticing that
%It is enough to prove that $\bsmat B^\top (-s\,I_n-A^\top)^{-1} & I_n \esmat
%L _{X}\,\bsmat (s\,I_n-A)^{-1}\,B \\[1mm] I_n \esmat=0$. Indeed
%\beann
%[ \begin{array}{cc} B^\top (-s\,I_n-A^\top)^{-1} & I_n   \end{array} ]\,
%L _{X}\,\bmat{c}  (s\,I_n-A)^{-1}\,B  \\[0mm] I_n \emat & = & -B^\top (s\,I_n+A^\top)^{-1}[(s\,I_n+A^\top)X  -  X(s\,I_n-A)](s\,I_n-A)^{-1} B   \\
%&& -B^\top (s\,I_n+A^\top)^{-1} X	\,B+B^\top X (s\,I_n-A)^{-1} B =0.
%\eeann
\beann
&& \hspace{-1cm} [ \begin{array}{cc} B^\top (-s\,I_n-A^\top)^{-1} & I_n   \end{array} ]
(\Pi _{X}-\Pi)\,\left[ \begin{array}{cc} (s\,I_n-A)^{-1}\,B \\ I_n \end{array} \right] \\
&& =-B^\top (s\,I_n+A^\top)^{-1}[(s\,I_n+A^\top)X-X(s\,I_n-A)](s\,I_n-A)^{-1} B \\
&& -B^\top (s\,I_n+A^\top)^{-1} X	\,B+B^\top X (s\,I_n-A)^{-1} B =0.
\eeann
\endproof

The following important result relates the rank of the spectrum $\Phi(s)$ with that of the matrix $R$, and it provides an explicit expression for a square spectral factor of $\Phi(s)$.

\begin{theorem}
Let $X=X^\top$ solve CGCARE($\Sigma$). Then, 
\begin{enumerate}
\item the normal rank of $\Phi(s)$ is equal to the rank of $R$;
\item $W(s) \defi R^{\frac{1}{2}} R^\dagger S _{X}^\top (s\,I_n-A)^{-1}B+R^{\frac{1}{2}}$ is a square spectral factor of $\Phi(s)$. 
\end{enumerate}
\end{theorem}
\proof
As already observed, since $X$ is a solution of CGCARE($\Sigma$), 
there holds $\ker R \subseteq \ker (XB)$. It follows that $\Pi _{X}$ can be written as 
$\Pi _{X}=V\,\bsmat Q _{X}-S _{X}\,R^\dagger S _{X}^\top  & 0 \\[1mm] 0 & R\esmat V^\top$, 
 where $V= \bsmat I_n & S _{X}\,R^\dagger \\[1mm] 0 & I_m \esmat$. Moreover, if $X$ solves CGCARE($\Sigma$), we get $Q _{X}-S _{X}\,R^\dagger S _{X}^\top=0$, and $\Pi _{X}$ can be factored as
$\Pi _{X}=\left( V \bsmat 0 & 0 \\[1mm] 0 & R^{\frac{1}{2}} \esmat \right)\left( \bsmat 0 & 0 \\[1mm] 0 & R^{\frac{1}{2}} \esmat V^\top\right)$. From $\bsmat 0 & 0 \\[1mm] 0 & R^{\frac{1}{2}} \esmat \bsmat I_n & 0 \\[1mm] R^\dagger\,S _{X}^\top & I_m \esmat=\bsmat 0 & 0 \\[1mm] R^{\frac{1}{2}}  R^\dagger\,S _{X}^\top & R^{\frac{1}{2}}  \esmat$, we find that $\Phi(s)$ can be written as $\Phi(s)=W^\top(-s)\,W(s)$, where
$W(s)=R^{\frac{1}{2}} R^\dagger S _{X}^\top (s\,I_n-A)^{-1}B+R^{\frac{1}{2}} = R^{\frac{1}{2}} [I_m+ R^\dagger S _{X}^\top (s\,I_n-A)^{-1}B]$.
Thus we can write $W(s)= R^{\frac{1}{2}}\,T _{X}(s)$, where $T _{X}(s) \defi I_m+ R^\dagger\,S _{X}^\top\,(s\,I_n-A)^{-1}B$ is square and invertible for all but finitely many $s \in \complex$. Its inverse can be written as 
$T^{-1} _{X}(s)=I_m-R^\dagger S _{X}^\top (s\,I_n-A _{X})^{-1}B$. Thus, the normal rank of $(T _{X}^{\top}(-s))^{-1} \Phi(s) T _{X}^{-1}(s)= R$
is equal to  the normal rank $r$  of $\Phi(s)$. %Let us define $B _{2}= B\,G$.
\endproof

In the following lemma, given a solution of CGCARE($\Sigma$), a subspace that will be shown to play a crucial role in the solution of the associated optimal control problem will be introduced. This subspace is the reachable subspace associated with the pair $(A _{X}, B\,G)$.

\begin{lemma}
Let $X=X^\top$ solve CGCARE($\Sigma$) and define
 \begin{equation}
\label{defR0X}
\gR _{0,X} \defi \ima [\begin{array}{ccccccccccc} B\,G && A _{X}\,B\,G && A _{X}^2\,B\,G&&  \ldots&&A _{X}^{n-1}\,B\,G \end{array}].
 \end{equation}
 %to be the reachable subspace associated with the pair $(A _{X}, B\,G)$. 
 Let $C _{X}\defi C-D\,R^\dagger\,S _{X}^\top$. There holds $\gR _{0,X} \subseteq \ker C _{X}$. 
\end{lemma}
\proof 
Since $\Phi(s)=G^\top(-s) G(s)=W^\top(-s)\,W(s)$ with $W(s)= R^{\frac{1}{2}}\,T _{X}(s)$, we find
\beann
G(s)\,T _{X}^{-1}(s) &=& \left(C\,(s\,I_n-A)^{-1}B+D\right) \left( I_m-R^\dagger S _{X}^\top (s\,I_n-A _{X})^{-1}B \right) \\
&=& C\,(s\,I_n-A)^{-1}B+D+C\,(s\,I_n-A _{X})^{-1}B \\
&&-D\,R^\dagger\,S _{X}^\top (s\,I_n-A _{X})^{-1}B \\
&=& (C-D\,R^\dagger\,S _{X}^\top) (s\,I_n-A _{X})^{-1}B+D,
\eeann
where the first equality follows from observing that $B\,R^\dagger S _{X}^\top=A-A _{X}=(s\,I_n-A _{X})-(s\,I_n-A)$. %The subscript $X$ has been dropped because this subspace is independent of the solution $X$ of CGCARE($\Sigma$).
We have already shown that $(T _{X}^{\top}(-s))^{-1} \Phi(s) T _{X}^{-1}(s)= R$. Thus, $\ker R \subseteq \ker G(s) T _{X}^{-1}(s)$. Hence,
$G(s) T _{X}^{-1}(s) \ker R=C _{X}\,(s\,I_n-A _{X})^{-1}B\ker R+D\,\ker R=\{0\}$. Since $D\,\ker R=\{0\}$, then $C _{X}\,(s\,I_n-A _{X})^{-1}B _{2}=0$. Therefore, $\gR _{0,X}$ must be in $\ker C _{X}$.
\endproof

In the case where  $X=X^\top$ is the solution of GCARE($\Sigma$) corresponding to the optimal cost, 
it is intuitive and simple to see that $\ker X$ is output-nulling for the quadruple $(A,B,C,D)$ and the corresponding gain $-K _{X}$ is a friend of $\ker X$, on the basis of the optimality and of the fact that the cost cannot be smaller than zero in view of the positivity of the index. 
Stated differently, if $x_0\in \ker X$, applying the control $u(t)=-K_X\,x(t)$ ensures that $x(t)\in \ker X$ for all $t \ge 0$, and the cost remains at zero, i.e.,
\[
\bmat{cc} A-B\,K_X \\ C-D\,K_X \emat\,\ker X \subseteq \ker X \oplus \{0\}.
\]
However, the following much stronger result holds.

\begin{theorem}
\label{th-inv+friend}
Let  $X=X^\top$ be a solution of GCARE($\Sigma$). 
Then, $\ker X$ is an output-nulling subspace of the quadruple $(A,B,C,D)$ and $-K _{X}$ is a friend of $\ker X$, or, equivalently, $\ker X$ is $A _{X}$-invariant and contained in the null-space of $C _{X}$.
\end{theorem}
\proof
{ Since $X$ is a solution of GCARE($\Sigma$), the closed-loop Lyapunov equation 
 \begin{equation}
\label{clgdare}
X\,A _{X}+A _{X}^\top X+Q _{0X}=0
 \end{equation}
holds, where 
$Q _{0X}   \defi    Q-S\,R^\dagger S^\top +X\,B\,R^\dagger B^\top X=  % (C^\top-S _{X} R^\dagger D^\top)(C-D\,R^\dagger S _{X}^\top) = 
C _{X}^\top C _{X} \ge 0$.
%\eeann
%$Q _{0X}\defi Q-S\,R^\dagger S^\top +X\,B\,R^\dagger B^\top X \geq 0$.} %In view of Lemma \ref{lemstein}, $\ker X$ is $A _{X}$-invariant and is contained in the null-space of $Q _{0X}$. 
%We show that $Q _{0X}=C _{X}^\top C _{X}$. Indeed,
%\beann
%C _{X}^\top\,C _{X} =& (C^\top-S _{X} R^\dagger D^\top)(C-D\,R^\dagger S _{X}^\top) = Q-S\,R^\dagger S _{X}^\top-S _{X} R^\dagger S^\top+S _{X} R^\dagger S _{X}^\top = Q _{0X}.% -S\,R^\dagger S^\top +X\,B\,R^\dagger B^\top X.
%\eeann
Moreover, from the definition of $C _{X}$ we also get
$Q _{0X}=C _{X}^\top\,C _{X}=[\,I_n \;\;\; -K _{X}\,]\Pi \bsmat I_n \\[1mm] -K _{X}^\top \esmat \geq 0$. Now, consider the Lyapunov equation $X\,A _{X}+A _{X}^\top X+C _{X}^\top\,C _{X}=0$, and let $\xi \in \ker X$. By multiplying this equation from the left by $\xi^\top$ and from the right by $\xi$ we obtain $C _{X}\,\xi=0$, which says that $\ker X \subseteq \ker C _{X}$. With this fact in mind, we multiply the same equation from the right by $\xi$, and we obtain $X\,A _{X}\,\xi=0$, which says that $\ker X$ is $A _{X}$-invariant. Thus, $\ker X$ is an $A _{X}$-invariant subspace contained in the null-space of $C _{X}$, and is therefore an output-nulling subspace for $(A,B,C,D)$ and $-K _{X}=-R^\dagger S _{X}$ is an associated  friend.
\endproof

We recall that we have defined the subspace $\gR _{0,X}$ as the reachability subspace of the pair $(A_X,B\,G)$. Since $A_X$ depends on the solutions $X=X^\top$ of CGCARE($\Sigma$) considered, at first glance it appears that the subspace $\gR _{0,X}$ also depends on $X$. However, we now prove that this is not the case:
 the subspace $\gR _{0,X}$
is independent of the particular solution $X=X^\top$ of CGCARE($\Sigma$). Moreover, $A _{X}$ restricted to this subspace does not depend on the particular
solution $X=X^\top$ of CGCARE($\Sigma$).

\begin{theorem}
\label{thee}
Let $X=X^\top$ be a solution of CGCARE($\Sigma$), and let $\gR _{0,X}$ be defined by (\ref{defR0X}). Then, 
\begin{itemize}
\item $X\,\gR _{0,X}=\{0\}$;
\item $\gR _{0,X}$ is independent of $X$;
\item $A _{X}|_{\gR _{0,X}}$ is independent of $X$.
\end{itemize}
\end{theorem}
\proof
Since $\gR _{0,X}$ is $A _{X}$-invariant and is contained in $\ker C _{X}$, in a basis of the state space adapted to $\gR _{0,X}$ we have
$\gR _{0,X}= \ima \bsmat  I_r   \\[1mm]   0    \esmat$, $A _{X}=\bsmat   A _{X,11}   & A _{X,12}   \\[1mm]   0   & A _{X,22}   \esmat$, $B _{2}=\bsmat  B _{21}   \\[1mm] 0 \esmat $, $C _{X}= [ \begin{array}{cccc} 0 & C _{X,2}\end{array} ]$,
where $r=\dim \gR _{0,X}$. If we partition $X$ conformably with this basis as $X=\bsmat X _{11} & X _{12} \\[1mm] X _{12}^\top & X _{22} \esmat$, we need to show that $X _{11}=0$ and $X _{12}=0$. 
Due to the structure of $C _{X}$, by pre- and post-multiplying the closed-loop Lyapunov equation $X\,A _{X}+A _{X}^\top X+C _{X}^\top\,C _{X}=0$ by $\bsmat I_r & 0 \esmat$ and $\bsmat I_r \\[1mm] 0 \esmat$, respectively,  we get
 $X _{11}\,A _{X,11}+A _{X,11}^\top X _{11}=0$. Now, $\ker R \subseteq \ker (X\,B)$ implies $X\,B _{2}=0$, which in turn implies  $X _{11}\,B _{21}=0$ and $X _{12}^\top B _{21}=0$. Therefore, $X _{11}$ satisfies
$\bsmat X _{11}\,A _{X,11}+A _{X,11}^\top X _{11} & X _{11} B _{21} \\[1mm] B _{21}^\top X _{11} & 0 \esmat=0$. Since the pair $(A _{X,11},B _{21})$ is reachable, it is always possible to choose a matrix $K$ in such a way that $A _{X,11}+B _{21}\,K$ has unmixed spectrum. Thus,
 \beann
 0& =&\bmat{cc}  I_r   &   K^\top   \\[-1mm]   0   &   I_{n-r}   \emat        \bmat{cc}   X _{11}\,A _{X,11}  +  A _{X,11}^\top X _{11} &  X _{11} B _{21}   \\[-1mm]   B _{21}^\top X _{11}   &   0   \emat      \bmat{cc}   I_r   &   0   \\[-1mm]   K   &   I_{n-r}   \emat \\
 & = &  \bmat{cc}   X _{11}(A _{X,11}  +  B _{21}K)  +  (A _{X,11}  +  B _{21}K)^\top X _{11}  &  X _{11} B _{21}   \\[-1mm]   B _{21}^\top X _{11}   &   0  \emat
 \eeann
 gives $X _{11}=0$. The Lyapunov equation (\ref{clgdare}) reads as
 \beann
 \bmat{cc} 0 & X _{12}\,A _{X,22} \\[-1mm] X _{12}^\top A _{X,11} & \star \emat+ \bmat{cc} 0 & A _{X,11}^\top X _{12} \\[-1mm] A _{X,22}^\top X _{12}^\top & \star \emat+
 \bmat{cc} 0 &0 \\[-1mm] 0 & \star \emat=0,
 \eeann
 which leads to
 $X _{12}\,A _{X,22}+A _{X,11}^\top X _{12}=0$. This identity, together with $X _{12}^\top B _{21}=0$, leads to $X _{12}=0$ in view of the observability of the pair $(A _{X,11}^\top,B _{21}^\top)$. Thus, $X\,\gR _{0,X}=\{0\}$. 
 
 We now want to show that $\gR _{0,X}$ is independent of $X$, where $X=X^\top$ is a solution of CGCARE($\Sigma$). %From (\ref{kercond}), we have already observed that it follows that $\ker R \subseteq \ker (X\,B)$. 
 In a certain basis of the input space, we can write $R=\bsmat R_{0} && 0 \\[1mm] 0 && 0 \esmat$, where $R_{0}$ is positive definite. Matrix $B$ can be written conformably with this basis as $B=[\,B_{1}\;\;B _{2}\,]$. From (\ref{kercond}), in this basis we must have $B _{2}^\top \,X=0$, i.e., $X\,B\,G=0$. Let us write $A _{X}=F-B\,R^\dagger\,B^\top X$, where $F\defi A-B\,R^\dagger\,S^\top$.
 We show that $\gR _{0,X}$, i.e, the reachable subspace of the pair $(A _{X},B\,G)$, coincides with that of the pair $(F,B\,G)$, which is independent of $X$ since $F$ is independent of $X$. First, we observe that
$A _{X}\,B _{2} =  (F-B\,R^\dagger\,B^\top X)\,B _{2}=F\,B _{2}$, since as already observed $X\,B _{2}=0$. We now prove by induction that $A _{X}^j\,B _{2}=F^j\,B _{2}$ for all $j \in \mathbb{N}$. The statement has been proved for $j=1$. Assume  $A _{X}^k\,B _{2}=F^k\,B _{2}$ for some $k >1$. 
First, in view of Theorem \ref{th-inv+friend}, $\ker X$ is $A _{X}$-invariant, which also implies that $A _{X}^k\,\ker X\subseteq \ker X$, i.e., $X\,A _{X}^k\,\ker X=\{0\}$. On the other hand, since $\ima B _{2}\subseteq \ker X$, we have also $X\,A _{X}^k\,B _{2}=X\,F^k\,B _{2}=0$. Thus,
$A _{X}^{k+1}\,B _{2}= A _{X}\,F^k\,B _{2}= (F-B\,R^\dagger\,B^\top X)\,F^k\,B _{2}=F\,F^k\,B _{2}=F^{k+1}\,B _{2}$.
 It is now clear that 
 \[
 \gR _{0,X}=\ima [ \begin{array}{ccccccccc}  B _{2} && A _{X}\,B _{2} &&  \ldots && A _{X}^{n-1}\,B _{2}\end{array} ]=
 \ima [ \begin{array}{ccccccccc}  B _{2} && F\,B _{2} && \ldots &&F^{n-1}\,B _{2}\end{array} ]
 \]
 which is independent of $X$. 
 We now prove that $A _{X}|_{\gR _{0,X}}$ is independent of $X$.
 Let $Y=Y^\top$ now be another solution of CGCARE($\Sigma$). Let $A _{Y}$ be the corresponding closed-loop matrix. We find
$A _{X}-A _{Y} = B\,R^\dagger (S _{Y}^\top-S _{X}^\top) =B\,R^\dagger\,B^\top (Y-X)$.
 We want to show that in this basis we have 
$A _{Y}=\bsmat A _{X,11} & A _{Y,12} \\[1mm] 0 & A _{Y,22} \esmat$.
 From the considerations above, since it has been already proved that
 $\gR _{0,X}=\gR_{0,Y}$,  in this basis we have $X=\bsmat 0 & 0 \\[1mm] 0 & X _{22} \esmat$ and $Y=\bsmat 0 & 0 \\[1mm] 0 & Y _{22} \esmat$, so that
 $A _{Y}=A _{X}-\bsmat \star \;&\; \star \\[1mm] \star  \;&\;  \star \esmat \bsmat 0 & 0 \\[1mm] 0 & X _{22}-Y _{22} \esmat=
 \bsmat A _{X,11} & \star\\[1mm] 0 & \star \esmat$,
 which shows that $A _{X}|_{\gR _{0,X}}=A _{Y}|_{\gR_{0,Y}}$. 
 \endproof

{The next result shows that the reachable subspace associated with the pair $(A _{X}, B\,G)$, which we denoted by $\gR _{0,X}$, coincides with the largest reachability output-nulling subspace on the output-nulling subspace $\ker X$. In view of Theorem \ref{thee}, such reachability output-nulling subspace (and the corresponding restriction of the closed-loop mapping to it) is therefore
independent of the particular solution $X=X^\top$ of CGCARE($\Sigma$) that we consider.}

\begin{theorem}
\label{the}
Let $X=X^\top$ be a solution of CGCARE($\Sigma$). Let $\gR^\star_{\ker X}$ be the largest reachability subspace on $\ker X$. Then,
$\gR^\star_{\ker X}=\gR _{0,X}$.
\end{theorem}
\proof
Since $\gR_{0}$ is the reachable subspace of the pair $(A _{X},B\,G)$, it is  the smallest $A _{X}$-invariant subspace containing $\ima (B\,G)=B\,\ker D$. On the other hand, the reachability subspace $\gR^\star_{\ker X}$ on $\ker X$ is  % {characterized as
%follows \cite[Theorem  7.14]{Trentelman-SH-01}, \cite[p.~424]{Ntogramatzidis-05}: Let $F$ be an {\em arbitrary} friend of $\ker X$, i.e., $F$ is any feedback matrix such that 
%$(A+B\,F)\ker X \subseteq  \ker X\subseteq \ker (C+D\,F)$. 
the smallest $(A+B\,F)$-invariant subspace 
containing $\ker X \cap B\,\ker D$, where $F$ is an {\em arbitrary} friend of $\ker X$, i.e., $F$ is any feedback matrix such that $(A+B\,F)\ker X \subseteq  \ker X\subseteq \ker (C+D\,F)$, \cite[Theorem  7.14]{Trentelman-SH-01}. 
 Notice that $\gR^\star_{\ker X}$ does not depend on the choice of the friend $F$, \cite[Theorem  7.18]{Trentelman-SH-01}.
We have seen in Theorem \ref{th-inv+friend} that $F= -K _{X}$ is a particular friend of $\ker X$.
For this choice of $F$, we have $A+B\,F=A-B\,K _{X}=A _{X}$, so that $\gR^\star_{\ker X}$ is the smallest $A _{X}$-invariant subspace 
containing $\ker X \cap B\,\ker D$. It is easy to see that $\ker X \cap B\,\ker D$ coincides with $B \ker D$, because $\ker X \supseteq B \ker D$ in view of the inclusion $\ker R \subseteq \ker XB$ following from (\ref{kercond}).
\endproof

{
\section{The Hamiltonian system}
\label{ESDE}
The aim of this section is to establish a link between the geometric properties of the solutions of CGCARE($\Sigma$) presented in the previous section and the structure of the so-called Hamiltonian system, which plays a crucial role in the study of the solutions of continuous-time (differential and algebraic) Riccati equations. %pencil $N-s\,M$, where
%Since in the present setting we are not assuming that $R$ is positive definite, the {\em trasversality} condition (\ref{el5}) cannot be solved in $u(t)$ so as to lead to a set of $2n$ equations in the sole variables $x(t)$ and $\lambda(t)$.
%However, a very convenient form in which equations (\ref{el1}), (\ref{el3}) and (\ref{el5}) can be written, that does not require inversion of $R$, is the  descriptor form
%\bea
%\label{eq6}
%M\,\dot{p}(t)=N\,p(t) \qquad t\,{\in}\,[0,T],
%\eea
%where
Recall that the Hamiltonian system associated with the Popov triple $\Sigma$ is defined by the equations
\bea
\begin{array}{rcl}
 \bmat{c} \dot{x}(t) \\[-1mm] \dot{\lambda}(t) \emat &=& \bmat{cc} A & 0 \\[-1mm] -Q & -A^\top \emat \bmat{c} x(t) \\[-1mm] \lambda(t)\emat+\bmat{c} B \\[-1mm] -S \emat u(t) \\
 y(t)&=&[\begin{array}{cc} S^\top & B^\top \end{array}] \bmat{c} x(t) \\[-1mm] \lambda(t)\emat+R\,u(t), \end{array} \label{ham}
\eea
where the variable $\lambda(t)$ is the costate vector. We define $\hat{A} \defi  \bsmat A & 0 \\[1mm] -Q & -A^\top \esmat$, $\hat{B} \defi  \bsmat B \\[1mm] -S \esmat$, $\hat{C} \defi \bsmat S^\top & B^\top \esmat$ and $\hat{D} \defi R$. The Hamiltonian system (\ref{ham}) is identified with the matrix quadruple $(\hat{A}, \hat{B}, \hat{C}, \hat{D})$. 
The Hamiltonian system has strong relations with the corresponding optimal control problem. Indeed, using an Euler-Lagrange approach, the optimality conditions of an LQ problem can be written as in (\ref{ham}) with the additional constraint $y(t) =0$ for all $t \ge 0$. %In the regular case (i.e., with $\det R \neq 0$),  the output equation $y(t)=0$ can be solved in $u(t)$, giving rize to an autonomous differential equation in the sole state-costate, and the corresponding dynamic matrix of this equation is referred to as the {\em Hamiltonian matrix}..
It is a classic and very well-known result that the set of invariant zeros %Equivalently, given the largest output-nulling subspace $\gV^\star$ of the quadruple $(A,B,C,D)$, i.e., the largest among all subspaces $\gV$ satisfying the inclusion $\bsmat A \\[1mm] C \esmat \gV \subseteq (\gV^\star \oplus \{0\}) +\ima \bsmat B \\[1mm] D \esmat$, the invariant zeros are those eigenvalues of $A+B\,F$ that do not depend on the particular $F$ for which $\bsmat A+B\,F \\[1mm] C+D\,F \esmat \gV^\star \subseteq \gV^\star \oplus \{0\}$. }
of the Hamiltonian system is mirrored with respect to the imaginary axis, see e.g. \cite{Prattichizzo-MN-04}. Moreover, given a solution $X$ of the standard continuous-time algebraic Riccati equation, the invariant zeros of the Hamiltonian system (\ref{ham}) are given by the union of the eigenvalues of the closed-loop matrix $A _{X}$ with those of $-A _{X}$. In symbols,
\begin{equation}
\label{std}
\gZ({\hat{A}, \hat{B}, \hat{C}, \hat{D}})= \sigma(A _{X}) \cup \sigma(-A _{X}).
\end{equation}
The goal of this section is to show that when $R$ is singular but the CGCARE($\Sigma$) admits a solution $X$, the set of invariant zeros of the Hamiltonian system (\ref{ham}) is a subset of such union.
More precisely, the following result holds.
\begin{theorem}
\label{thefund}
Let $X$ be a solution of CGCARE($\Sigma$). Let the pair $(A _{X},B\,G)$ be written in the reachability form as
$\bsmat 
A _{X,11} & A _{X,12} \\[1mm] 0 & A _{X,22} \esmat$, $\bsmat B_{2,1} \\[1mm] 0 \esmat$, 
where the pair $(A _{X,11},B_{2,1})$ is completely reachable. Let $\Gamma _{X} \defi A _{X,22}$.
%
%
%Let $\Gamma _{X}$ represent the uncontrollable eigenvalues of the pair $(A _{X},B\,G)$, where $G$ is a basis for $\ker R$. 
There holds
\[
\gZ(\hat{A}, \hat{B}, \hat{C}, \hat{D})=\sigma(\Gamma _{X}) \cup \sigma(-\Gamma _{X}).
\]
\end{theorem}
\ \\[-5mm]
In order to prove Theorem \ref{thefund}, we need the following technical lemmas.
\ \\[-5mm]
\begin{lemma}
\label{zeri}
The set of invariant zeros of a quadruple $(A,B,C,D)$ is invariant with respect to state feedback and output injection and with respect to changes of coordinates in the state space, i.e., for any matrices $F$ and $G$ and for any non-singular $T$ of suitable sizes
 there hold
\beann
\gZ(A,B,C,D)&=&\gZ(A+B\,F,B,C+D\,F,D)\\
&=&\gZ(A+G\,C,B+G\,D,C,D) \\
&=&\gZ(T^{-1}A\,T,T^{-1} B, C\,T,D).
\eeann
%Moreover, the set of invariant zeros of a quadruple $(A,B,C,D)$ is invariant with respect to changes of coordinates in the state space, i.e., for any non-singular $T$ of suitable size
%\[
%\gZ(A,B,C,D)=\gZ(T^{-1}A\,T,T^{-1} B, C\,T,D).
%\]
\end{lemma}
\proof
The first equality follows by observing that for all matrices $F$ and for all $s \in \complex$ there holds
$\bsmat A+B\,F-s\,I_n & B \\[1mm] C+D\,F & D \esmat = \bsmat  A-s\,I_n & B \\[1mm] C & D \esmat \bsmat  I_n & 0 \\[1mm] F & I_m \esmat$. %, 
%$\bsmat  A+G\,C-s\,I_n & B+G\,D \\[1mm] C & D \esmat= \bsmat  I_n & G \\[1mm] 0 & I_p \esmat\bsmat  A-s\,I_n & B \\[1mm] C & D \esmat$.
The second is dual. The third statement follows from
$\bsmat T^{-1}A\,T-s\,I_n & T^{-1}B \\[1mm] C\,T & D \esmat = \bsmat  T^{-1} & 0 \\[1mm] 0 & I_p \esmat\bsmat  A-s\,I & B \\[1mm] C & D \esmat \bsmat  T & 0 \\[1mm] 0 & I_m \esmat$.
\endproof
\begin{lemma}
\label{lemnuovo}
Let  $X=X^\top$ be a solution of CGCARE($\Sigma$). 
The invariant zeros of the Hamiltonian system (\ref{ham}) coincide with the generalized eigenvalues of the matrix pencil
\be\label{tiblokdec}\hat{P}(s)=\left[ \begin{array}{ccc}
A _{X}-s\,I_n & 0 & B \\[-1mm]
0 & -(A _{X}^\top+s\,I_n) & 0 \\[-1mm]
0 & B^\top &  R \end{array} \right].
\ee
\end{lemma}
%Let $\Gamma _{X}$ represent the uncontrollable eigenvalues of the pair $(A _{X},B\,G)$, where $G$ is a basis for $\ker R$. Therefore,
%\[
%\gZ(\hat{\mathfrak{Q}})=\sigma(\Gamma _{X}) \cup \sigma(-\Gamma _{X}).
%\]
%
\proof 
We perform a state-feedback transformation in (\ref{ham}).
Let $u(t)=[\begin{array}{cc} -K _{X} & 0 \end{array} ] \bsmat x(t) \\[1mm] \lambda(t) \esmat+v(t)$, so that
\beann
&& \bmat{c} \dot{x}(t) \\[-1mm]  \dot{\lambda}(t) \emat = \bmat{cc} A-B\,K _{X} & 0 \\[-1mm]  -Q+S\,K _{X} & -A^\top \emat \bmat{c} x(t) \\[-1mm]  \lambda(t)\emat+\bmat{c} B \\[-1mm]  -S \emat v(t) \\
&& y(t)=[\begin{array}{cc} S^\top-R\,K _{X}  & B^\top \end{array} ]\bmat{c} x(t) \\[-1mm]  \lambda(t)\emat+R\,v(t)
\eeann
Now we change coordinates in the state-space of the Hamiltonian system with $T=\bsmat I_n & 0 \\[1mm] X & I_n
\esmat$, and we obtain
\beann
 \hat{A}^\prime & = &  T^{-1} \hat{A}\,T= \bmat{cc} I_n & 0 \\[-1mm]  -X & I_n
\emat \bmat{cc} A-B\,K _{X} & 0 \\[-1mm]  -Q+S\,K _{X} & -A^\top \emat \bmat{cc} I_n & 0 \\[-1mm]  X & I_n
\emat \\
&=& \diag \{A _{X},-A^\top\},
\eeann
since $-X\,A _{X}-Q+S\,K _{X}-A^\top X=0$ in view of CGCARE($\Sigma$). Moreover
$\hat{B}^\prime=T^{-1}\,\hat{B}=\bsmat I_n & 0 \\[1mm] -X & I_n
\esmat \bsmat B \\[1mm] -S \esmat=\bsmat B \\[1mm] -X\,B-S \esmat$
and
$\hat{C}^\prime = \hat{C}\,T= [\begin{array}{cc} S^\top-R\,K _{X}  & B^\top\end{array}]  \bsmat I_n & 0 \\[1mm] X & I_n
\esmat =
[\begin{array}{ccc} S^\top-R\,R^\dagger S _{X}  && B^\top\end{array}]
=[ \begin{array}{ccc} 0  && B^\top \end{array} ]$.
%because $X$ solves the {\bf C}-GCARE. Therefore, the matrices of the Hamiltonian system become
%\beann
%&& A^\prime=\bmat{cc} A _{X} & 0 \\ 0 & -A^\top \emat \quad B^\prime=\bmat{c} B \\ -X\,B-S \emat \\
%&& C^\prime = \bmat{cc} 0  & B^\top \emat \quad D^\prime =R
%\eeann
Finally, $\hat{D}^\prime =R$. In view of Lemma \ref{zeri}, we have $\gZ(\hat{A}, \hat{B}, \hat{C}, \hat{D})=\gZ(\hat{A}^\prime,   \hat{B}^\prime, \hat{C}^\prime, \hat{D}^\prime)$.
Now we perform an output-injection using the matrix $G=\bsmat 0 \\[1mm] K _{X}^\top \esmat$
and we obtain 
\beann
 \hat{A}^\second & = &  \hat{A}^\prime +G\,\hat{C}=\bmat{cc}  A _{X}   & 0   \\[-1mm]   0   & -A^\top   \emat +\bmat{c} 0 \\    K _{X}^\top    \emat [ \begin{array}{cc} 0  & B^\top \end{array} ]
=
\bmat{cc} A _{X} & 0 \\[-1mm]  0 & -A _{X}^\top \emat\\
 \hat{B}^\second & = &  \hat{B}^\prime +G\,\hat{D}=\bmat{c} B \\[-1mm]  -X\,B-S \emat+\bmat{c} 0 \\[-1mm]  K _{X}^\top \emat R=
\bmat{c} B \\[-1mm]  0 \emat \\
 \hat{C}^\second & = & \hat{C}^\prime= [ \begin{array}{ccc} 0  && B^\top \end{array} ] \quad  \qquad \hat{D}^\second = \hat{D}^\prime=R,
\eeann
where we have used the fact that $-X\,B-S+S _{X}\,R^\dagger R=0$ since $S\,R^\dagger R=S$ and $XB(R^\dagger R-I_m)=0$.
Again, in view of Lemma \ref{zeri}, we have $\gZ(\hat{A}, \hat{B}, \hat{C}, \hat{D})=\gZ(\hat{A}^\second,   \hat{B}^\second, \hat{C}^\second, \hat{D}^\second)$. Thus, the invariant zeros of the Hamiltonian system $(\hat{A}, \hat{B}, \hat{C}, \hat{D})$ are the values of $\zeta \in \complex$ such that the Rosenbrock matrix pencil (\ref{tiblokdec}) loses rank.
\endproof
%\hat{P}(s)=\left[ \begin{array}{ccc}
%A _{X}-s\,I_n & 0 & B \\
%0 & -(A _{X}^\top+s\,I_n) & 0 \\
%0 & B^\top &  R \end{array} \right].
%\ee
%loses rank with respect to its normal rank.
%
%

It is worth remarking that the generalized eigenvalues of $\hat{P}(s)$ are independent of the solution $X=X^\top$ of CGCARE($\Sigma$), since these coincide with the invariant zeros of the Hamiltonian system.
Observe also that when $R$ is non-singular {(i.e., when $X$ is a solution of CARE($\Sigma$))}, 
this result allows us to re-obtain (\ref{std}), since clearly 
 \begin{equation}
\label{division}
\sigma(\hat{P}(s))=\sigma(A _{X}-s\,I_n) \cup \sigma (A _{X}^\top+s\,I_n),
 \end{equation}
 where the symbol $\sigma(\hat{P}(s))$ stands for the set of generalized eigenvalues of the pencil $\hat{P}(s)$ counting multiplicities. However, (\ref{division}) does not hold when $R$ is singular.
}
 \begin{example}
 \label{ex0}
 {
Let
$A=\bsmat -4 & 0 \\[1mm] 2 & 6 \esmat$, $B=\bsmat 0 & -7 \\[1mm] 2 & -4 \esmat$, $Q=\bsmat  \frac{17}{4} & 0 \\ 0 & 0  \esmat$, $S=\bsmat 0 & 0 \\[1mm] 0 & 0 \esmat$, $R=\bsmat  0 & 0 \\[1mm] 0 & 4  \esmat$.
The matrix $X=\diag\{-1,0\}$ is a solution of CGCARE($\Sigma$) but  CARE($\Sigma$) is not defined in this case.  The closed-loop matrix is $A _{X}=\bsmat 33/4 & 0 \\[1mm] 9 & 6 \esmat$. Applying the result in Lemma \ref{lemnuovo} we find that the Rosenbrock matrix associated with the Hamiltonian system can be written as % \\[-11mm]
{\small
\beann
\hat{P}(s)=\left[ \begin{array}{cc|cc|cc}
 \frac{33}{4} - s &  0 &     0 &  0 & \;0 & -7 \\[-1mm]
     9 & 6-s &     0 &  0 & \;2 & -4 \\
     \hline
     0 &  0 & -\frac{33}{4}-s &  -9 & \;0 & 0\\[-1mm]
     0 &  0 &     0 & -6-s & \;0 & 0\\[-1mm]
     \hline
     0 &  0 &   0 &  2 & \;0 & 0\\[-1mm]
     0 &  0 &     -7 &  -4 &\; 0 & 4 \end{array} \right].
     \eeann
     }
The normal rank of $\hat{P}(s)$ is equal to $5$. The eigenvalues of $A _{X}$ are equal to $33/4$ and $6$. While it is true that when $s=\pm 33/4$ the rank of $\hat{P}(s)$ is equal to $4$, for both $s=6$ and $s=-6$ the rank of $\hat{P}(s)$ is equal to $5$.
This result says that, unlike the regular case, not all the eigenvalues of $A _{X}$ are invariant zeros of the Hamiltonian system. Specifically, the invariant zeros of the Hamiltonian system are $\pm 33/4$.
% \footnote{{We warn that the routine {\tt eig.m} of the software MATLAB$^{\textrm{\tiny{\textregistered}}}$ (version 7.11.0.584(R2010b)) in this case fails to provide the right answer when used  for the computation of the generalised eigenvalues of the pencil $N-s\,M$. It indeed returns $1$ as a generalized eigenvalue of the pencil.} }
 \endex
 }
 \end{example}

%\begin{itemize}
%\item $A _{X}$ is singular if and only if $R$ or $A-B\,R^\dagger\,S$ is singular. Indeed, from (\ref{newN}) it is clear that the singular part of $N$ is divided between the singular part of $A _{X}$ and that of $R _{X}$, i.e.,
%\[
%\ker N = A _{X}^{-1}\,B\,\ker R _{X} \oplus 0_n \oplus \ker R _{X},
%\]
%where $A _{X}^{-1}$ denotes the inverse of the map $A _{X}$. Therefore, the fact that $N$ is singular does not necessarily imply that $A _{X}$ is singular. Moreover, since the singular part of $N$ does not depend on the solution $X=X^\top$ of CGDARE($\Sigma$), and neither does that of $R _{X}$, then also the singular part of $A _{X}$ is independent of the solution $X=X^\top$ of CGDARE($\Sigma$);
% the generalized eigenvalues of the extended symplectic matrix pencil are given by
%\beann
%\sigma(A _{X}-s\,I_n) \cup \sigma \left( \left[ \begin{array}{cc} I_n-s\,A _{X}^\top & 0 \\
% -s\,B^\top &  R _{X}\end{array} \right] \right) 
% \eeann
% \end{itemize}

\begin{theorem}
\label{th0}
Let $X$ be a solution of CGCARE($\Sigma$). Two matrices ${U} _{X}$ and ${V} _{X}$ exist such that
 %\bea
%&&\hspace{-11cm} {U} _{X}\,\hat{P}(s)\,{V} _{X} = \nonumber 
%\eea
%\ \\[-20mm]
{\small
\bea
 && \hspace{-3mm}{{\Large {U} _{X}\,\hat{P}(s)\,{V} _{X} =}} \nonumber \\
 && \left[ \begin{array}{cc|c|ccc}
 A _{X,11}-s\,I_{r} & B _{21} & 0 & A _{X,12} & 0 & B _{11} \\
 \hline
 0 & 0 & -(A _{X,11}^\top+s\,I_r) & 0 & 0 & 0 \\[-1mm]
 0 & 0 & B _{21}^\top & 0 & 0 & 0 \\
 \hline
 0 & 0 & 0 & A _{X,22}-s\,I_{n-r}\!\!& 0 & B _{12} \\[-1mm]
 0 & 0 & -A _{X,12}^\top & 0 &  \!\! -(A _{X,22}^\top+s\,I_{n-r})   & 0 \\[-1mm]
 0 & 0 & B _{11}^\top & 0 & B _{12}^\top & R_{0} 
  \end{array} \right], \label{fc}
  \eea
  }
  where the pair $(A _{X,11},B _{21})$ is reachable and $R_{0}$ is invertible. Moreover, the sub-matrix pencil 
  \[
  \hat{P}_{1}(s) \defi \bsmat
 A _{X,22}-s\,I_{n-r}& 0 & B _{12} \\
 0 & -(A _{X,22}^\top+s\,I_{n-r}) & 0 \\
  0 & B _{12}^\top & R_{0} 
\esmat
\]
 in (\ref{fc})
is regular, and the generalized eigenvalues of the pencil $\hat{P}(s)$ are the generalized eigenvalues of $\hat{P}_{1}(s)$.%\\[-5mm]
  \end{theorem}
\proof
Consider an orthogonal change of coordinate in the input space $\real^m$ induced by {the $m \times m$ orthogonal matrix $T=[\,T_{1}\;\;T _{2}\,]$ where $\ima T_{1}=\ima R$ and $\ima T_{2}=\ima G=\ker R$.}  In this basis $R$ is block-diagonal, with the first block being non-singular and the second being zero, i.e.,
$T^\top R\,T= \diag \{R_{0},0\}$,
 where $R_{0}$ is invertible. Its dimension is denoted by $m_{1}$. Consider the block matrix $\hat{T}\defi\diag\{I_n, I_n, T\}$. % Let $U _{X}$ and $V _{X}$ be as defined in the proof of Theorem \ref{the0}. 
 By multiplying $\hat{P}(s)$ on the left by $\hat{T}^\top$ and on the right by $\hat{T}$, and by defining  the matrices $B_{1} \defi B\,T_{1}$ and $B _{2} \defi B\,T_{2}$ we get%\\[-15mm]
 {\small
 \beann
 \hat{T}^\top\hat{P}(s)\, \hat{T}= 
 \left[ \begin{array}{cccc}
 A _{X}-s\,I_n & 0 & B_{1} & B _{2} \\[-1mm]
0 & -(A _{X}^\top+s\,I_n) & 0 & 0\\[-1mm]
0 & B_{1}^\top &  R_{0} & 0 \\[-1mm]
0 & B _{2}^\top &  0 & 0 
 \end{array} \right].
 \eeann
 }
 Notice that $\ima B _{2}= \ima (B\,G)$ in view of the identity $\ker R=\ima G$. Matrix $B_{1}$ has $m_{1}$ columns. Let us denote by $m _{2} \defi m - m_{1}$ the number of columns of $B _{2}$. Let us now take a matrix
 % At this point we perform a Kalman reachability canonical decomposition on the pairs $(A _{X},B _{2})$ and $(A _{X}^\top,B _{2}^\top)$ by taking a matrix 
 $H=[\,H_{1}\;\;H _{2}\,]$ such that $\ima H_{1}$ is the reachable subspace from the origin of the pair $(A _{X},B _{2})$, which coincides with the subspace $\ker R$, yielding
$H^{-1}\,A _{X}\,H=\bsmat   A _{X,11}   &   A _{X,12}    \\
   0   &   A _{X,22}  \esmat   , \; H^{-1}\,B _{2}=\bsmat    B _{21}   \\[1mm]   0    \esmat$, $H^{-1}\,B_{1}=\bsmat    B _{11}   \\[1mm]   B _{12}   \esmat$.
  Let $\hat{H}=\diag\{H,H,I_{m_{1}},I_{m _{2}}\}$ be partitioned conformably with the block structure of the pencil. %We obtain
% \beann
% &&\hat{H}^{-1} \hat{T}^\top\hat{P}(s)\, \hat{T} \hat{H}=  \nonumber 
 %\eeann
% {\small
% \beann
%&&  \left[ \begin{array}{cc|cc|cc}
% A _{X,11}-s\,I_{r} & A _{X,12} & 0 & 0 & B _{11} & B _{21} \\
% 0 & A _{X,22}-s\,I_{n-r} & 0 & 0 & B _{12} & 0\\
% \hline
%0 & 0 & -(A _{X,11}^\top+s\,I_r) & 0 & 0 & 0\\
%0 & 0 & -A _{X,12}^\top & -(A _{X,22}^\top+s\,I_{n-r}) & 0 & 0 \\
%\hline
%0 & 0 & B _{11}^\top & B _{12}^\top & R_{0} & 0 \\
%0 & 0 & B _{21}^\top & 0 &  0 & 0 
% \end{array} \right],
 %\eeann
 %}
 %where $r$ is the size of the reachable subspace of the pair $(A _{X},B _{2})$. 
Reordering $\hat{H}^{-1} \hat{T}^\top\hat{P}(s)\, \hat{T} \hat{H}$ via two suitable unimodular matrices $\Omega_{1}$ and $\Omega_{2}$
% \beann
% \Omega_{1} \defi \left[ \begin{array}{cccccc}
 %  I_r   &  0   &  0   &  0   &  0   &  0   \\
%    0   &  0   &  I_r   &  0   &  0   &  0   \\
 %    0   &  0   &  0   &  0   &  0   &  I_{m _{2}}   \\
%      0   &  I_{n-r}   &  0   &  0   &  0   &  0   \\
%       0   &  0   &  0   &  I_{n-r}   &  0   &  0   \\
%        0   &  0   &  0   &  0   &  I_{m_{1}}   &  0 
%      \end{array} \right] \qquad \textrm{and} \qquad
 %      \Omega _{2} \defi \left[ \begin{array}{cccccc}
 %  I_r   &  0   &  0   &  0   &  0   &  0   \\
 %   0   &  0   &  0   &  I_{n-r}   &  0   &  0   \\
 %    0   &  0   &  I_r   &  0   &  0   &  0   \\
%      0   &  0   &  0   &  0   &  I_{n-r}   &  0   \\
%       0   &  0   &  0   &  0   &  0   &  I_{m_{1}}   \\
%        0   &  I_{m _{2}}   &  0   &  0   &  0   &  0 
  %     \end{array} \right],
  %     \eeann
       yields (\ref{fc}) with $\hat{U} _{X}=\Omega_{1}\,\hat{H}^{-1}\,\hat{T}^\top$ and $\hat{V} _{X} \hat{T}\,\hat{H}\, \Omega _{2}$, where $r$ is the size of the reachable subspace of the pair $(A _{X},B _{2})$.
 { We now proceed with the computation of the normal rank of $P(s)$ .}   
%\begin{lemma}
  %\label{lemf}
 %Consider the partitioned matrix
 %\beann
 %\Xi=\left[ \begin{array}{cc} \Xi _{11} & \Xi _{12} \\ O & \Xi _{22} \end{array} \right]
 %\eeann
%{If either $\Xi _{11}$ is full row-rank or  $\Xi _{22}$ is full column-rank, then $\rank \Xi=\rank \Xi _{11}+\rank \Xi _{22}$.}
 %\end{lemma}
Since the pair $(A _{X,11}, B _{21})$ is reachable by construction, all the $r$ rows of the submatrix $[\,A _{X,11}-s\,I_{r} \;\; B _{21} \,]$ are linearly independent for every $s \in \complex \cup \{\infty\}$. This also means that of the $r+m _{2}$ columns of $[\,A _{X,11}-s\,I_{r} \;\; B _{21} \,]$, only $r$ are linearly independent, and this gives rise to the presence of a null-space of $P(s)$ whose dimension $m _{2}$ is independent of $s \in \complex \cup \{\infty\}$. Thus,
\beann
\rank \hat{P}(s)=r+
\rank \left[ \begin{array}{c|ccc}
  -(A _{X,11}^\top+s\,I_r) & 0 & 0 & 0 \\[-1mm]
 B _{21}^\top & 0 & 0 & 0 \\
 \hline
0 & A _{X,22}-s\,I_{n-r}& 0 & B _{12} \\[-1mm]
 -A _{X,12}^\top & 0 & -(A _{X,22}^\top+s\,I_{n-r}) & 0 \\[-1mm]
  B _{11}^\top & 0 & B _{12}^\top & R_{0} 
  \end{array} \right].
  \eeann
  Again, since the pair $(A _{X,11}, B _{21})$ is reachable, then $(A _{X,11}^\top, B _{21}^\top)$ is observable, and the rank of the submatrix $\left[ \begin{smallmatrix} 
 -(A _{X,11}^\top+s\,I_r)  \\ B _{21}^\top \end{smallmatrix} \right]$ is constant and equal to $r$ for every $s \in \complex \cup \{\infty\}$. Thus,
$\rank \hat{P}(s)=2\,r+ \rank P_{1}(s)$, 
  where %\\[-11mm]
  {\small
  \beann
   \hat{P}_{1}(s) = \left[ \begin{array}{ccc}A _{X,22}-s\,I_{n-r}& 0 & B _{12} \\[-1mm]
 0 & -(A _{X,22}^\top+s\,I_{n-r}) & 0 \\[-1mm]
  0 & B _{12}^\top & R_{0} 
  \end{array} \right].
  \eeann
  }
 Since $\det \hat{P}_{1}(s)=\det(A _{X,22}-s\,I_{n-r}) \cdot \det(-(A _{X,22}^\top+s\,I_{n-r}))\cdot  \det R_{0}$, a value $s \in \complex$ can certainly be found for which $\det \hat{P}_{1}(s) \neq 0$. This means that the normal rank of $\hat{P}_{1}(s)$ is equal to $2\,(n-r)+m_{1}$, and therefore
$\normrank \hat{P}(s)=2\,r+2\,(n-r)+m_{1}=2\,n+m_{1}$.
 It also follows that the generalized eigenvalues of the pencil $\hat{P}(s)$ are the values $s \in \complex \cup \{\infty\}$ for which the rank of $\hat{P}_{1}(s)$ is smaller than its normal rank $2\,(n-r)+m_{1}$. These values are the eigenvalues of $A _{X,22}$ plus their opposites, including possibly the eigenvalue at infinity, whose multiplicity {| be it algebraic or geometric | is %, in general, {\em not} given by the sum of $m_{1}$ plus the multiplicity of  the eigenvalue in zero of $A _{X,22}$.
%In fact, the multiplicity of the eigenvalue at infinity is 
the 
multiplicity of the zero eigenvalue of the matrix}
$P_{\infty}\defi \diag\{I_{n-r},I_{n-r},0_{m_{1}}\}$.
  {The last $m_{1}$ columns of $P_{\infty}$ give rise to an eigenvalue at infinity whose multiplicity (algebraic and geometric) is exactly equal to $m_{1}$, since in this case the dimension of $\ker P_{\infty}$ is equal to $m_{1}$.   
\endproof
  
 Theorem \ref{thefund} now follows as a corollary of Theorem \ref{the}. Indeed, from (\ref{fc}) we find
$\gZ(\hat{A}, \hat{B}, \hat{C}, \hat{D})=\sigma(A _{X,22}) \cup \sigma(-A _{X,22})$.
It turns out that, unlike the regular case, not all the eigenvalues of the closed-loop matrix $A _{X}$ are invariant zeros of the Hamiltonian system (\ref{ham}). In particular, the eigenvalues of $A _{X}$ restricted to $\gR^\star_{\ker X}$ -- which are the controllable eigenvalues of the pair $(A _{X},B\,G)$ -- are not invariant zeros of the Hamiltonian system, whereas the eigenvalues induced by $A _{X}$ on  $\real^n /  \gR^\star_{\ker X}$ along with their opposites are invariant zeros of the Hamiltonian system.

 \begin{example}
 {
 Consider Example \ref{ex0}. Using the solution $X=\diag \{-1,0\}$ of CGCARE($\Sigma$) we easily find that $\ker R$ and $\ima R$ are respectively spanned by the vectors $\bsmat 1 \\[1mm] 0 \esmat$ and $\bsmat 0 \\[1mm] 1 \esmat$.  Hence, by taking $T=\left[ \begin{smallmatrix} 0 & 1 \\[1mm] 1 & 0 \end{smallmatrix} \right]$ we obtain
$T^\top R\,T=\diag\{ 4,0\}$.
 Hence, in this case $m_{0}=1$. Moreover, we partition $B\,T$ as $B\,T=\left[ \begin{smallmatrix} -7 & 0 \\[1mm] -4 & 2 \end{smallmatrix} \right]$, so that $B_{1}=\left[ \begin{smallmatrix} -7  \\[1mm] -4  \end{smallmatrix} \right]$ and $B _{2}=\left[ \begin{smallmatrix} 0 \\[1mm] 2 \end{smallmatrix} \right]$.
As expected, the image of $B _{2}=B\,G$ is exactly equal to the reachability subspace on $\ker X=\spanR \left\{\left[ \begin{smallmatrix} 0 \\[1mm] 1 \end{smallmatrix} \right]\right\}$, which in this case coincides with $\spanR \left\{\left[ \begin{smallmatrix} 0 \\[1mm] 1 \end{smallmatrix} \right]\right\}$. The normal rank of the $\hat{P}(s)$ is equal to $2\,n+m_{0}=5$.
The invariant zeros of the Hamiltonian system are given by the uncontrollable eigenvalues of the pair $(A _{X},B _{2})=\left(\left[ \begin{smallmatrix} 33/4 & 0 \\[1mm] 9 & 6 \end{smallmatrix} \right],\left[ \begin{smallmatrix} 0 \\[1mm] 2 \end{smallmatrix} \right]\right)$ plus their opposites, i.e., $33/4$ and one at $-33/4$. Since $A _{X,22}=0$ and $B _{12}=2$, the matrix pencil $\hat{P}(s)$ also has a generalized eigenvalue at infinity with multiplicity equal to the dimension of $\ker \left[ \begin{smallmatrix} 0 & 0 \\[1mm] 2 & 0 \end{smallmatrix} \right]$, which is equal to $1$. %Finally, the multiplicity of the normal null-space of the extended symplectic pencil is equal to $1$.
By writing the Rosenbrock matrix pencil associated with the Hamiltonian system in the form given by (\ref{fc}), we get in fact %\\[-11mm]
{\small
 \bea
 \hat{T}^\top \hat{P}(s)\, \hat{T}= 
 \left[ \begin{array}{cc|c|ccc}
  6-s & 2 & 0 & 9 & 0 & -4 \\
  \hline
    0 & 0 & -6-s & 0 & 0 & 0 \\[-1mm]
      0 & 0 & 2 & 0 & 0 & 0 \\
      \hline
        0 & 0 & 0 & \frac{33}{4}-s & 0 & -7 \\[-1mm]
          0 & 0 & -9 & 0 & -(\frac{33}{4}+s) & 0 \\[-1mm]
            0 & 0 & -4 & 0 & -7 & 4 \end{array} \right], \label{tgy}
 \eea
 }
 which  shows that $33/4$ and $-33/4$ are indeed the only finite generalized eigenvalues of $\hat{P}(s)$.
% from which it is clear that zero is a generalized eigenvalue.
% Notice also that the generalized eigenvalue of the extended symplectic pencil is equal to the eigenvalue of $A _{X}$ restriced to the quotient space $\ker X / \gR^\star_{\ker X}$ plus the eigenvalues of $A _{X}$ restricted to the quotient space $\real^n / \ker X$. In this case, since $\gR^\star_{\ker X}=\ker X$, there are no eigenvalues of $A _{X}$ restriced to $\ker X / \gR^\star_{\ker X}$. The only eigenvalue of $A _{X}$ restricted to the quotient space $\real^n / \ker X$ is exactly equal to zero. In fact, let us take a change of basis $T=[\,T_{1}\;\;\;T _{2}\,]$ such that $\ima T_{1}=\ker X$. For example, let us take $T_{1}$ equal to the first canonical basis vector of $\real^2$ and $T _{2}$ equal to the second canonical basis vector. Hence, $T^{-1}\,A _{X}\,T=A _{X}$, and the entry at the bottom right corner is exactly the single eigenvalue of $A _{X}$ restricted to the quotient space $\real^n / \ker X$. \endex
 }
 \end{example}

 \begin{remark}
  {
  The
MATLAB$^{\textrm{\tiny{\textregistered}}}$ routine for the solution of the continuous-time algebraic Riccati equation is {\tt care.m}. This routine requires matrix $R$ to be positive definite, and delivers the stabilizing solution of this equation (which exists if and only if $(A,B)$ is stabilizable and the Hamiltonian matrix has no eigenvalues on the imaginary axis). Thus, {\tt care.m} cannot handle the case where $R$ is singular. Using the decomposition in this section, and using {\tt care.m} for the regular part of the Hamiltonian pencil delivers the solution of CGCARE($\Sigma$) which, loosely speaking, is {\em as stabilizing as possible}. Differently from what happens in the standard case, when no stabilizing solutions of CGCARE($\Sigma$) exist, it may still be possible to add another feedback which stabilizes the system, as the following section will show.
}
\end{remark}

{
\section{Stabilization}
 In the previous sections, we have observed that the eigenvalues of the closed-loop matrix $A _{X}$ restricted to the subspace $\gR_{0,X}$ are independent of the particular solution $X=X^\top$ of CGCARE($\Sigma$) considered. This means that these eigenvalues -- which do not appear as invariant zeros of the Hamiltonian system -- are present in the closed-loop regardless of the solution $X=X^\top$ of CGCARE($\Sigma$) that we consider.
 On the other hand, we have also observed that $\gR_{0,X}$ coincides with the subspace $\gR^\star_{\ker X}$, which is by definition the smallest $(A-B\,K _{X})$-invariant subspace containing $\ker X \cap B\,\ker D=\ima (B\,G)$. It follows that it is always possible to find a matrix $L$ that assigns all the eigenvalues of the map $(A _{X}+B\,G\,L)$ restricted to the reachable subspace $\gR^\star_{\ker X}$, by adding a further term $B\,G\,L\,x(t)$ to the feedback control law, because this does not change the value of the cost with respect to the one obtained by $u(t)=-K _{X}\,x(t)$. Indeed, the additional term only affects the part of the trajectory on $\gR^\star_{\ker X}$ which is output-nulling. However, in doing so it may stabilize the closed-loop if $\ker X$ is externally stabilized by $-K _{X}$, see \cite{Trentelman-SH-01}. Indeed, since
 $\gR_{0,X}$ is output-nulling with respect to the quadruple $(A,B,C,D)$, it is also output-nulling for the quadruple $(A-B\,K _{X},B,C-D\,K _{X},D)$, and two matrices $\Xi$ and $\Omega$ exist such that
 \begin{equation}
 \label{XU}
 \bmat{c} A _{X} \\[-1mm] C _{X} \emat \,R_{0,X}=\bmat{c} R_{0,X} \\[-1mm]  0 \emat \Xi+\bmat{c} B \\[-1mm]  D \emat \Omega, 
 \end{equation}
 where $R_{0,X}$ is a basis matrix of $\gR_{0,X}$.  In order to find a feedback matrix that stabilizes the system, we solve the former in $\Xi$ and $\Omega$, so as to find $L$ such that
 $\bsmat A _{X}+B\,L \\[1mm] C _{X}+D\,L \esmat \,R_{0,X}=\bsmat R_{0,X} \\[1mm] 0 \esmat \Xi$,
  where the eigenvalues of $\Xi$ are the eigenvalues of the map $A _{X}+B\,L$ restricted to $\gR_{0,X}$.
  %Using the standard procedure of geometric control theory \cite{Wonham-85,Basile-M-92,Trentelman-SH-01}, w
  We first compute the set of solutions of (\ref{XU}) in $\Xi$ and $\Omega$, i.e.,
  \bea
  \label{hu}
  \bmat{c} \Xi \\[-1mm] \Omega \emat=\bmat{c} \hat{\Xi} \\[-1mm] \hat{\Omega} \emat+\bmat{c} H_{1}\\[-1mm] H _{2} \emat K,
  \quad
  \bmat{c} \hat{\Xi} \\[-1mm] \hat{\Omega} \emat\defi \bmat{cc} R_{0} & B \\[-1mm] 0& D \emat^\dagger \bmat{c} A _{X} \\[-1mm] C _{X} \emat R_{0}
  \eea
  and $\bsmat H_{1} \\[1mm] H _{2} \esmat$ is a basis matrix of $\ker \bsmat R_{0} & B \\[1mm] 0 & D \esmat$. Since $\gR_{0,X}$ is a controllability subspace, the pair $(\hat{\Xi},H_{1})$ is reachable. This implies that a matrix $K$ in (\ref{hu}) can always be found so that the eigenvalues of $\Xi$ are freely assignable (provided they come in complex conjugate pairs). Hence, we use such $K$ in (\ref{hu}) and then we compute $L=-\Omega\,R_{0}^\dagger$. This choice guarantees that only the eigenvalues of $A _{X}$ restricted to $\gR_{0,X}$ get affected by the use of $L$.
   %
 %We show this fact in the following example.
%

 \begin{example}
 \label{ex2}
  {
  Let $A = \bsmat -8  &  0\\[1mm]
     6   & 0 \esmat$, $B= \bsmat 0 \\[1mm]
     -4 \esmat$, $C = [\begin{array}{cc} 4    &               0  \end{array} ]$, $D=0$, so that
     $Q=C^\top C=\diag \{16,0\}$, $S=C^\top D=0$ and $R=D^\top D=0$.    % \textcolor{red}{Let $X=\bsmat x_{1} & x _{2} \\ x _{2} & x_3 \esmat$. Then
    % \[
    % X\,A+A^\top\,X-(S+X\,B)\,R^{\dagger}\,(S^\top+B^\top X)+Q=0,
    % \]
    % becomes
    %\[
    % \bmat{cc} x_{1} & x _{2} \\ x _{2} & x_3 \emat \bmat{cc} -8  &  0\\
    % 6   & 0 \emat + \bmat{cc} -8  &  6\\
    % 0   & 0 \emat \bmat{cc} x_{1} & x _{2} \\ x _{2} & x_3 \emat+\bmat{cc} 16 & 0 \\ 0 & 0 \emat=0
    % \]
    % because $R^\dagger=0^\dagger=0$. This leads to
    % \beann
    % && -8\,x_{1}+6\,x _{2}-8\,x_{1}+6\,x _{2}+16=0 \\
    % && -8\,x _{2}+6\,x_3=0
    % \eeann
    % i.e., 
    % \beann
    % && x _{2}=\frac{3}{4}\,x_3 \\
    % && x_{1}=\frac{9}{16}\,x_3+1 
    % \eeann     
    % }
          %
     %
    One can directly verify that the set of solutions of GCARE($\Sigma$) is parameterized in $t \in \real$ as
    $X_t= \bsmat \frac{9}{16}\,t+1 & \frac{3}{4}\,t \\[1mm] \frac{3}{4}\,t & t \esmat$. Thus, $X=X_{0}=\bsmat 1 & 0 \\[1mm] 0 & 0 \esmat$ is the only solution of GCARE($\Sigma$) for which $\ker R \subseteq \ker (S+X\,B)$. This implies that $X=X_{0}=\bsmat 1 & 0 \\[1mm] 0 & 0 \esmat$ is the only solution of CGCARE($\Sigma$), and it is positive semidefinite.
    % for any $t \in \real$, the matrix $X_t= \bsmat \frac{9}{16}\,t+1 & \frac{3}{4}\,t \\[1mm] \frac{3}{4}\,t & t \esmat$
  %   is a solution of CGCARE($\Sigma$), since $X$ satisfies (\ref{gcare}) and $\ker R=\ker (S+X\,B)$. Notice that $X=X_{0}=\bsmat 1 & 0 \\[1mm] 0 & 0 \esmat$ is the smallest positive semidefinite solution of CGCARE($\Sigma$). 
  %   In this case $G=I_m -R^\dagger R=1$, which means that $B _{2}=B$.
     Since $R=0$, we find $K _{X}=[\,0\;\;\;0\,]$, which implies $A _{X}=A$. Hence, CGCARE($\Sigma$) does not admit a stabilizing solution. However, we now see that the infinite-horizon problem admits an optimal  solution which is also stabilizing. Indeed, we find $\gR_{0,X}=\ima \bsmat 0 \\[1mm] 1 \esmat$. The eigenvalue of $A _{X}$ restricted to $\gR_{0,X}$ is 0, while the eigenvalue induced by the map $A_{X}$ on the quotient space $\real^2/\gR_{0,X}$ is $-8$. 
     The optimal trajectory is
     \[
     \bmat{c} x_{1}(t) \\[-1mm] x _{2}(t)\emat=e^{A _{X}\,t}\bmat{c} x_{1}(0) \\[-1mm] x _{2}(0)\emat=\bmat{cc} e^{-8\,t} & 0 \\[-1mm] \frac{3}{4}\left( 1-e^{-8\,t}\right) & 1 \emat \bmat{c} x_{1}(0) \\[-1mm] x _{2}(0)\emat,
     \]
     which implies that the optimal cost is $J^*=x_{1}^2(0)$,
 %    \beann
 % J^*=\int_{0}^\infty 16\,e^{-16\,t}\,x_{1}(t)\,dt=x_{1}^2(0),
  %\eeann  
  % \int_{0}^\infty  x^\top(0) \bmat{cc} e^{-8\,t} & \frac{3}{4}\left( 1-e^{-8\,t}\right) \\[-1mm] 0 & 1 \emat  Q \bmat{cc} e^{-8\,t} & 0 \\[-1mm] \frac{3}{4}\left( 1-e^{-8\,t}\right) & 1 \emat  x(0)  dt, %\\
%     & = & \bmat{cc} x_{1}(0) & x _{2}(0) \emat   \bmat{cc} 16\,e^{-16\,t} & 0 \\ 0 & 0 \emat_{0}^{\infty} \bmat{c} x_{1}(0) \\ x _{2}(0) \emat=x_{1}^2(0),
   %  \eeann
  i.e., it coincides with $x^\top(0)\,X_{0}\,x(0)$. We can find another optimal solution that assigns the additional eigenvalue of the closed loop to $-1$. In this case, $\bsmat R_{0,X} & B \\[1mm] 0 & D \esmat^\dagger=\bsmat 0 & 0 \\ 1 & -4 \\ 0 & 0 \esmat^\dagger=\frac{1}{17} \bsmat 0 & 1 & 0 \\[1mm] 0 & -4 & 0 \esmat$, so that $\bsmat \hat{\Xi} \\[1mm] \hat{\Omega}\esmat=\bsmat 0 \\[1mm] 0 \esmat$ using (\ref{hu}). Moreover, a basis for the null-space of $\bsmat R_{0,X} & B \\[1mm] 0 & D \esmat$ is $\bsmat 4 \\[1mm] 1 \esmat$. We find 
$\bsmat \Xi \\[1mm] \Omega \esmat=\bsmat 4 \\[1mm] 1 \esmat \,K$.      Imposing $\Xi=-1$ gives $K=-1/4$, which in turn gives 
$L=-\Omega\,R_{0,X}^\dagger=-K\,\bsmat 0 \\[1mm] 1 \esmat^\dagger=\frac{1}{4} [\begin{array}{cc} 0 & 1 \end{array}]=[\begin{array}{cc} 0 & \frac{1}{4} \end{array}]$.
     %using $L=[\,0\;\;\;1/4\,]$. We get $A+B\,L=\bsmat -8 & 0 \\[1mm] 6 & -1 \esmat$, so that
     Thus, $e^{(A+B\,L)\,t}=\bsmat e^{-8\,t} & 0 \\[1mm] \star & e^{-t} \esmat$, and the value of the cost remains $J^*=x_{1}^2(0)$.
   %  \beann
   %  J^*& = & \int_{0}^\infty [\begin{array}{cc} x_{1}(0) & x _{2}(0) \end{array}]   \bmat{cc}  e^{-8\,t}  & \star  \\  0  & e^{-t}  \emat   Q   \bmat{cc}  e^{-8\,t}  & 0  \\  \star  & e^{-t}  \emat     \bmat{c}   x_{1}(0)   \\   x _{2}(0)   \emat dt \\
    % & = & x_{1}^2(0).
    % \eeann
     This solution is optimal, and is also stabilizing. Thus, we found a stabilizing optimal control even in a situation in which CGCARE($\Sigma$) does not admit a stabilizing solution.
     }
     \end{example}

 \begin{remark}
  {
  The same procedure used in Example \ref{ex2} can be used also in examples where the eigenvalues of $A_{ X}$ are complex. Consider e.g. $A = \bsmat 1  &  1 \\[1mm]
     -1 & 1    \esmat$, $B= \bsmat 1 \\[1mm]
   0 \esmat$, and $Q$, $S$ and $R$ are zero matrices. The only solution of CGCARE($\Sigma$) is $X = 0$, so that $\sigma(A_{X})=\sigma(A)=\{1\pm i\}$. However, using the same procedure of Example \ref{ex2} we can find a matrix $L=[\begin{array}{cc} -9 & 19 \end{array}]$ which stabilizes the system since  $\sigma(A_{X}+B\,G\,L)=\{-3,-4\}$. Thus, 
   an optimal feedback that is stabilizing exists. 
   }
     \end{remark}

\begin{remark}
The case discussed in the previous example is somehow extreme. In fact, if CGCARE($\Sigma$) admits a solution and $R=0$ (which clearly implies $S=0$) it is clear that $BG=B$ and, for any solution $X$ 
of CGCARE($\Sigma$), $A_X=A$. Therefore in this case there exists a matrix $L$ such that the system can be stabilized by the feedback $A_X+BGL$ (that does not affect the cost index)  if and only if $(A,B)$ is stabilizable. 
This is an extreme case, but, as shown in the following example,  $R=0$ is far from being a necessary condition for the occurrence of cases where CGCARE($\Sigma$) admits solutions, none of which is stabilizing, but there exist a solution $X$ and a matrix $L$ such that $A_X+BGL$ is a stability matrix.   
\end{remark}

 \begin{example}
 \label{ex6.7}
  {
  Let $A = \bsmat 1  &  1 & 1\\[1mm]
     -3 & 1    & 0 \\[1mm] 1 & 0 & 0 \esmat$, $B= \bsmat 0 & 2 \\[1mm]
     0 & 0\\[1mm] 1 & 0 \esmat$, $Q = \bsmat 1  &  0 & -1\\[1mm]
     0 & 0    & 0 \\[1mm] -1 & 0 & 1 \esmat$, $S=\bsmat 1 & 0 \\[1mm] 0 & 0 \\[1mm] -1 & 0 \esmat$ and $R=\bsmat 1 & 0 \\[1mm] 0 & 0 \esmat$. 
      One can directly verify that the only two solutions of CGCARE($\Sigma$) are
    $X_{0}=0$ and $X_{1}= \diag\{0,0,2\}$. 
        None of these two solutions is stabilizing. Indeed, the eigenvalues of the closed-loop matrix relative to $X_{0}$ are $\{1,1\pm i\,\sqrt{3}\}$, while those of the one relative to $X_{1}$ are $\{-1,1\pm i\,\sqrt{3}\}$. Thus, CGCARE($\Sigma$) does not have a stabilizing solution.
Let us consider the solution  $X=X_{1}$. We have %, we obtain
$$
A_{X_1}=\bmat{cc|c}1\;&\;  1 \; &\; 1 \\ -3 \;&\;  1 \;&\;  0 \\ \hline 0 \;&\;  0 \;&\;  -1 \emat \qquad \text{and} \qquad BG=\bmat{c|c} 0 \;&\; 2 \\ 0 \;&\;  0 \\ \hline 0 \;&\;  0 \emat.
$$
Thus, by suitably selecting $L$ we can arbitrarily place the eigenvalues of the north-east corner of 
$ A_{X_1} +BGL$ while the third eigenvalue is fixed (and stable) and this new feedback does not affect the cost
(\ref{costinf}). For example, we can take
  \[
  L=\bmat{ccc} 0 & 0 & 0 \\  -\frac{7}{2} & \frac{3}{2} & 0\emat
  \]
  so that the overall closed-loop matrix becomes
  \[
  A_{X_1} +B\,G\,L=\bmat{ccc} -6 \;&\;  4 \;&\;  1 \\ -3 \;&\;  1 \;&\;  0 \\ 0 \;&\;  0 \;&\;  -1 \emat,
  \]
  whose eigenvalues are $\{-1,-2,-3\}$ and, hence, it is stable.     }
     \end{example}

\section*{Concluding remarks and future directions}
In this paper we investigated some structural properties of CGCARE arising from  singular LQ optimal control problems. These considerations revealed that a subspace can be identified that is independent of the particular solution of the Riccati equation considered and, even more importantly, such that the closed-loop matrix restricted to this subspace does not depend on the particular solution of the Riccati equation. %, and has been shown to be fixed for any state-feedback control constructed from a solution of the CGDARE. On the other hand, 
If such subspace is not zero, in the optimal control a further term can be added to the state feedback associated with the solution of the Riccati equation that does not affect the value of the cost function. This term can be expressed in state-feedback form, and can be used as a degree of freedom to be employed to stabilize the closed-loop even in cases in which no stabilizing solutions exist of the Riccati equation.
As for the discrete-time case, see \cite{Ferrante-04,NF-SCL-arxiv}, our analysis is expected to lead to a procedure for the order reduction of the CGCARE, which we believe will provide a relevant numerical   
edge in the solutions of CGCARE.

 %The solutions of (\ref{reduced}) are linked to the solutions of the extended symplectic difference equation
 %\beann
 %\bmat{ccc} I_n & O & O \\
 %O & -Z^\top & O \\
 %O & -B_2^\top & O \emat \bmat{c} \Xi_{t+1} \\ \Lambda_{t+1} \\ \Omega_{t+1} \emat =
 %\bmat{ccc} Z & O & B_2 \\
 %O & -I_n & O \\
 %O & O & R_0 \emat \bmat{c} \Xi_{t} \\ \Lambda_{t} \\ \Omega_{t} \emat.
 %\eeann
 %On the other hand, this equation can be written a suitable basis of the input space that isolates the invertible part of $R_0$ from its nilpotent part. In other words, the former can be written as
 %\beann
 %\bmat{cccc} I_n & O & O & O \\
 %O & -Z^\top & O & O \\
 %O & -B_{21}^\top & O & O \\
 %O & -B_{22}^\top & O & O
 %\emat \bmat{c} \Xi_{t+1} \\ \Lambda_{t+1} \\ \Omega_{t+1}^\prime \\ \Omega_{t+1}^\second \emat =
 %\bmat{cccc} Z & O & B_{21} & B_{22} \\
 %O & -I_n & O & O \\
 %O & O & \bar{R}_{0} & O \\
 %O & O & O & O \emat \bmat{c} \Xi_{t} \\ \Lambda_{t} \\ \Omega_{t}^\prime \\ \Omega_t^\second \emat.
 %\eeann

%


\begin{thebibliography}{100}

\bibitem{Abou-Kandil-FIJ-03}
H.~Abou-Kandil, G.~Freiling, V.~Ionescu and G.~Jank.
\newblock {\em Matrix {R}iccati Equations in Control and Systems Theory}.
\newblock Birkh{\"a}user, Basel, 2003.

%\bibitem{Colaneri-F-SCL}
%P.~Colaneri and A.~Ferrante.
%\newblock {A}lgebraic {R}iccati equation and {J}-spectral factorization in
%  ${H}_{\infty}$ estimation.
%\newblock {\em Systems \& Control Letters}, 41(5):383--393, April 2004.

\bibitem{Damm-04}
T. Damm.
\newblock {\em Rational Matrix Equations in Stochastic Control}.
\newblock Lecture Notes in Control and Information Sciences, Vol. 297, Springer Verlag, Berlin, 2004.


\bibitem{Damm-H-01}
T. Damm and D. Hinrichsen.
\newblock Newton's method for a rational matrix equation occurring in stochastic control.
\newblock {\em Linear Algebra Appl.}, vol. 332--334, pp. 81--109, 2001.

\bibitem{Dragan-MS-10}
V. Dragan, T. Morozan and A.-M. Stoica
\newblock {\em Mathematical Methods in Robust Control of Discrete-Time Linear Stochastic Systems}.
\newblock Springer, New York, 2010. 

\bibitem{Ferrante-04}
A.~Ferrante.
\newblock On the structure of the solution of discrete-time algebraic {R}iccati
  equation with singular closed-loop matrix.
\newblock {\em IEEE Transactions on Automatic Control}, AC-49(11):2049--2054,
  2004. 

\bibitem{NF-SCL-arxiv}
L.~Ntogramatzidis and A. Ferrante,
\newblock ``The discrete-time generalized algebraic Riccati
equation: order reduction and solutions' structure''.
\newblock {\em Systems \& Control Letters}, vol. 75, pp. 84--93, 2015. 

\bibitem{Ferrante-N-12}
A.~Ferrante, and L.~Ntogramatzidis,
\newblock ``The generalized discrete algebraic Riccati equation in linear-quadratic optimal
control''.
\newblock {\em Automatica}, vol. 49, pp. 471--478, 2013. 
%

%
%\bibitem{Ferrante-N-13-1}
%A.~Ferrante, and L.~Ntogramatzidis,
%\newblock ``A reduction technique for discrete generalized algebraic and difference Riccati equations'',
%\newblock {\em Linear and Multilinear Algebra}, vol. 62, pp. 1460--1474, 2014.

  \bibitem{Ferrante-N-14}
A.~Ferrante, and L.~Ntogramatzidis,
\newblock ``The generalized continuous algebraic Riccati equation and impulse-free continuous-time LQ optimal control''.
\newblock {\em Automatica}, vol. 50, pp. 1176--1180, 2014.

  \bibitem{Ferrante-N-14-1}
A.~Ferrante, and L.~Ntogramatzidis,
\newblock ``Continuous-time singular Linear-Quadratic control: necessary and sufficient conditions for the existence of regular solutions''.
\newblock {\em Systems \& Control Letters}, 93:30--34, 2016.

  \bibitem{Freiling-H-03}
G.~Freiling, and A.~Hochhaus,
\newblock ``Properties of the solutions of rational matrix difference equations''.
\newblock {\em Comput. Math. Applic.}, vol. 45, pp. 1137--1154, 2003.


  \bibitem{Freiling-H-04}
G.~Freiling, and A.~Hochhaus,
\newblock ``On a class of rational matrix differential equations arising in stochastic control''.
\newblock {\em Linear Algebra Appl.}, vol. 379, pp. 43--68, 2004.


\bibitem{Hautus-S-83}
M.L.J. Hautus and L.M. Silverman.
\newblock System structure and singular control.
\newblock {\em Linear Algebra Appl.}, vol. 50, pp. 369--402, 1983.

\bibitem{Kalaimani-BC-13}
R.K.~Kalaimani, M.N.~Belur and D.~Chakraborty
\newblock ``Singular LQ Control, Optimal PD Controller and Inadmissible Initial Conditions''.
\newblock {\em IEEE Trans. Aut. Control}, vol. 58, pp. 2603--2608, 2013.


\bibitem{Lancaster-95}
P.~Lancaster and L.~Rodman.
\newblock {\em Algebraic {R}iccati equations}.
\newblock Clarendon Press, Oxford, 1995.



\bibitem{Ionescu-O-96-1}
V.~Ionescu and C.~Oar\v{a}.
\newblock Generalized continuous-time {R}iccati theory.
\newblock {\em Linear Algebra Appl.}, vol. 232, pp. 111--130, 1996.


\bibitem{Ionescu-OW-99}
V.~Ionescu, C.~Oar\v{a}, and M.~Weiss.
\newblock {\em Generalized {R}iccati theory and robust control, a {P}opov
  function approach}.
\newblock Wiley, 1999.



%
\bibitem{Prattichizzo-MN-04}
D.~Prattichizzo, L.~Ntogramatzidis, and G.~Marro,
\newblock ``A new approach to the cheap LQ regulator exploiting the geometric properties of the Hamiltonian system''.
\newblock {\em Automatica}, vol. 44, pp. 2834--2839, 2008.
%

\bibitem{Rami-CZ-02}
A.M. Rami, X. Chen, and X.Y. Zhou, 
\newblock ``Discrete-time indefinite LQ control with state and control dependent noises''.
 \newblock {\em Journal of Global Optimization}, vol. 23, pp:~245--265, 2002.
 

\bibitem{Saberi-S-87}
A.\ Saberi and P.~Sannuti.
\newblock Cheap and singular controls for linear quadratic regulators.
\newblock {\em IEEE Trans. Aut. Control}, vol. 32, pp. 208--219,
 1987.


\bibitem{Saberi-SC-95}
A.~Saberi, P.~Sannuti, and B.M. Chen.
\newblock {\em $H_2$ Optimal Control}.
\newblock System and Control Engineering. Prentice Hall International, London,
  1995.
  
  
  \bibitem{Stoorvogel-92_2}
A.A. Stoorvogel.
\newblock The singular ${H}_2$ control problem.
\newblock {\em Automatica}, 28(3):627--631, 1992.
%

\bibitem{Stoorvogel-S-98}
A.A. Stoorvogel and A.~Saberi.
\newblock The discrete-time algebraic {R}iccati equation and linear matrix
  inequality.
\newblock {\em Linear Algebra Appl.}, vol. 274, pp. 317--365, 1998.


  %
\bibitem{Trentelman-SH-01}
H.L. Trentelman, A.A. Stoorvogel, and M.~Hautus.
\newblock {\em Control theory for linear systems}.
\newblock Springer,  2001.



\bibitem{Willems-KS-86}
J.C.~Willems, A.~K\`{\i}tap\c{c}i, and L.M.~Silverman.
\newblock ``Singular optimal control: a geometric approach".
\newblock {\em SIAM J. Control Optim.}, vol. 24, pp. 323--337, 
  1986.

\bibitem{Zhou-DG-96}
K.~Zhou, J.~Doyle, and K.~Glover.
\newblock {\em Robust and Optimal Control}.
\newblock Prentice Hall, New York, 1996.


\bibitem{zorzi1} M. Zorzi.
\newblock ``Convergence analysis of a family of robust Kalman filters based on the contraction principle".
\newblock {\em SIAM J. Control Optim.} To appear. 2017.

\bibitem{zorzi2} B. C. Levy, M. Zorzi. 
\newblock ``A contraction analysis of the convergence of risk-sensitive filters". 
\newblock {\em SIAM J. Control Optim.}, 54(4), 2154--2173, 2016.

\bibitem{zorzi3} M. Zorzi, B. C. Levy.
\newblock ``On the Convergence of a Risk Sensitive like Filter". 
{\em 54th IEEE Conference on Decision and Control}, Osaka, Japan, Dec. 2015. 


\end{thebibliography}
\end{document}